\documentclass [a4paper,twoside,13pts]{article}
\usepackage [francais]{babel}
\usepackage [latin1] {inputenc}
\usepackage [T1] {fontenc}
\usepackage{amssymb}
\usepackage{amsmath}
\usepackage{stmaryrd}
\usepackage{graphicx}
\title{Pénalisations de l'araignée brownienne (Penalizations of Walsh
Brownian motion) }
\author {Joseph Najnudel}
\begin{document} 
\maketitle 
\noindent

\textbf{Résumé : }
 Dans cet article, nous pénalisons la loi d'une araignée brownienne $(A_t)_{t
 \geq 0}$ prenant ses valeurs dans un ensemble fini $E$ de demi-droites concourantes, avec un poids égal
 à $\frac{1}{Z_t} \exp (\alpha_{N_t}
 X_t + \gamma L_t)$, où $t$ est un réel positif, $(\alpha_k)_{k \in
 E}$ une famille de réels indexés par $E$, $\gamma$ un paramètre réel, $X_
t$ la distance de $A_t$ à l'origine, $N_t$ ($\in E$) la demi-droite sur
 laquelle se trouve
 $A_t$, $L_t$ le temps local de $(X_s)_{0 \leq s \leq t}$ à l'origine,
 et $Z_t$ la constante de normalisation. Nous montrons que la famille des mesures
 de probabilité obtenue par ces pénalisations converge vers une probabilité
 limite quand $t$ tend vers l'infini, et nous étudions quelques propriétés
 de cette probabilité limite. \\ \\
\textbf{Abstract : } In this paper, we penalize a Walsh Brownian
 motion $(A_t)_{t \geq 0}$ (also called Brownian spider), which takes
 values in a finite set $E$ of intersecting rays, with a weight
 equal to $\frac{1}{Z_t} \exp (\alpha_{N_t} X_t + \gamma L_t)$, where $t$ is a
 positive real, $(\alpha_k)_{k \in E}$ a family of real numbers
 indexed by $E$, $\gamma$ a real parameter, $X_t$ the distance from
 $A_t$ to the origin, $N_t$ ($\in E$) the ray on which $A_t$ is to be found, $X_t$ the
 local time of $(A_s)_{0 \leq s \leq t}$ at the origin, and $Z_t$ the
 normalization constant. We show that the family of the probability measures
 obtained by these penalizations converges to a limit probability measure as $t$
 tends to infinity, and
 we study some properties of this limit probability
 measure. 
\\ \\
\textbf{Mots-clé : } pénalisation, temps local, araignée
 brownienne. 
\textbf{classifications AMS : } 60B10, 60J65 (60G17, 60G44, 60J25,
 60J55).  
\section{Présentation du problème et des principaux résultats obtenus}
\subsection{Introduction}
\noindent
Récemment, de nombreuses études de pénalisations du
mouvement brownien ont été effectuées, en particulier par B. Roynette, P. Vallois et
M. Yor (voir \cite{7}, \cite{8}, \cite{9}, \cite{10}, \cite{11}). \\
Dans \cite{9}, les pénalisations étudiées sont des fonctions de la valeur $X_t$
atteinte par un mouvement brownien en un temps $t$, et de
$S_t$, suprémum sur $[0,t]$ de ce mouvement
brownien. Plus précisément, on considère une
famille de mesures de probabilité $(\mathbf{W}^{(t)})_{t \geq 0}$
sur $\mathcal{C} (\mathbf{R}_+, \mathbf{R})$ vérifiant, pour tout
$\Lambda_t$ appartenant à la tribu $\mathcal{F}_t$ engendrée
par $(X_s)_{s \in [0,t]}$ ($(X_t)_{t \geq 0}$ étant le processus
canonique de $\mathcal{C} (\mathbf{R}_+, \mathbf{R})$) :
$$\mathbf{W}^{(t)}(\Lambda_t) = \frac{\mathbf{W} [\mathbf{1}_{\Lambda_t}
  f(X_t,S_t)]}{\mathbf{W} [f(X_t,S_t)]} $$ où 
$S_t$ est le maximum de $X_s$ pour $s \in [0,t]$, $\mathbf{W}$ la
mesure de Wiener, et $f$ une fonction de $\mathbf{R}^2$ dans
$\mathbf{R}_+$. \\ 
B. Roynette, P. Vallois et M. Yor montrent alors que pour certains
  choix de la fonction $f$, il existe une mesure de probabilité
  $\mathbf{W}^{(\infty)}$ (dépendant de $f$) sur $\mathcal{C} (\mathbf{R}_+, \mathbf{R})$
  telle que pour tout $s \geq 0$ et tout $\Lambda_s \in \mathcal{F}_s$
  : $$ \mathbf{W}^{(t)} (\Lambda_s) \underset{t \rightarrow
  \infty}{\rightarrow} \mathbf{W}^{(\infty)} (\Lambda_s)$$
Un des cas où cette convergence a lieu est celui où $f(a,y)= \exp(\lambda y + \mu a)$
  avec $\lambda$, $\mu \in \mathbf{R}$. \\ \\  
Par un changement de mouvement brownien, les résultats de \cite{9}
peuvent être adaptés au cas où $S_t$ est remplacé par $L_t$ (temps
local en 0 de $(X_u)_{u \leq t}$), et $X_t$ par $L_t -
|X_t|$; en effet, le théorème d'équivalence de Lévy affirme que $(S_t -
X_t, S_t)_{t \geq 0}$ a même loi que $(|X_t|,L_t)_{t \geq 0}$. \\ \\
Dans ces conditions, les poids exponentiels étudiés dans \cite{9} prennent
la forme : $\frac{1}{Z_t} \exp (\alpha |X_t| + \gamma L_t)$ où $\alpha$ et $\gamma$
sont des paramètres réels, et $Z_t$ est la constante de
normalisation. \\ \\
Le but de notre article est de généraliser l'étude de ces
pénalisations exponentielles à toutes les araignées browniennes
  prenant leurs valeurs dans un ensemble fini de demi-droites
concourantes.
\subsection{Quelques rappels et définitions}
\noindent
Dans ce paragraphe, nous allons définir le cadre général dans lequel on peut construire
  les araignées browniennes (étudiées dans \cite{1} et \cite{12}), et nous énoncerons
  plusieurs propriétés de ces processus, utiles par la suite. \\ \\
a) Soit $(E,\mu)$ un espace de probabilité fini; on suppose $\mu(\{ m \} ) >
0$ pour tout $m \in E$. Cet espace de probabilité est fixé une fois
  pour toutes dans cet article; par conséquent, nous
  omettrons en général d'indiquer la dépendance en $(E,\mu)$ des
quantités et des mesures de probabilités que nous introduirons. \\
On considère, sur l'espace $\mathbf{R}_E = \{ (0,0) \} \cup
(\mathbf{R}_+^* \times E)$, la distance $d$ définie par :
$$d((x,k),(y,l)) = |x-y| \mathbf{1}_{k = l} + (x+y) \mathbf{1}_{k \neq l}$$  
Cette distance permet de considérer $\mathcal{C}_E$, espace des
fonctions continues de $\mathbf{R}_+$ dans  $\mathbf{R}_E$, et de
munir cet espace de la tribu $\mathcal{T}_E$ associée à la topologie
de la convergence uniforme. \\ \\
b) Nous désignons par $(A_t = (X_t,N_t))_{t \geq 0}$ le processus canonique (à
valeurs dans $\mathbf{R}_E$) associé à l'espace $(\mathcal{C}_E,
\mathcal{T}_E)$ et nous notons, pour tout $t \in \mathbf{R}_+$,
$\mathcal{F}_t$ la sous-tribu de $\mathcal{T}_E$ engendrée par
$(A_s)_{0 \leq s \leq t}$. \\ 
Pour $(x,k) \in \mathbf{R}_E$, on peut alors considérer, sur
$\mathcal{C}_E$, la mesure de probabilité $\mathbf{W}_{(x,k)}$, sous
laquelle $(A_t)_{t \geq 0}$ est une araignée brownienne issue de 
$(x,k)$. \\ \\
c) Rappelons (voir \cite{1}) que cette araignée brownienne est un processus
de Feller qu'il est possible de caractériser par son semi-groupe
$(P_t)_{t \geq 0}$; pour toute fonction $f$ borélienne
bornée : $$ P_t f(x,k) = 2 \underset{m \in E}{\sum} \mu_m
\underset{\mathbf{R}_+^*}{\int} dy p_t(x+y) f(y,m) +
\underset{\mathbf{R}_+^*}{\int} dy (p_t(x-y) - p_t(x+y)) f(y,k)$$ 
avec $\mu_m = \mu( \{ m \} )$ (notation conservée dans la suite de l'article) et $p_t(a) = \frac{1}{\sqrt{2 \pi t}}
e^{-a^2/2t}$. \\ \\
d) Pour tout $(x,k) \in \mathbf{R}_E$, le processus  $(X_t)_{t \geq
  0}$, sous  $\mathbf{W}_{(x,k)}$, est un
mouvement  brownien réfléchi issu de $x$. \\  
D'autre part, si $T_0 = \inf \{ t \geq 0, X_t=0 \}$ et si
$\mathcal{I}$ est l'ensemble des intervalles d'excursion de $(X_t)_{t
  \geq T_0}$, $N_t$ est constant sur chaque intervalle $I \in
\mathcal{I}$ et on peut donc poser $N_t = N_I$ pour $t \in I$. On
montre alors que conditionnellement à $(X_t)_{t \geq 0}$, les
$(N_I)_{I \in \mathcal{I}}$ sont des variables aléatoires
indépendantes de loi $\mu$. \\
En particulier, pour tout $t \geq 0$, conditionnellement au fait que
$T_0 \leq t$, $N_t$ est une variable aléatoire de loi $\mu$, indépendante de $(X_s)_{s
  \geq 0}$. \\ \\ 
e) Dans notre étude de l'araignée brownienne, interviennent des processus à valeurs
réelles appelés processus bang-bang. \\
Par définition, un processus bang-bang de paramètre $\gamma > 0$ est
un processus $(Y_t)_{t \geq
  0}$, supposé issu de zéro dans cet article, et vérifiant l'équation différentielle stochastique : $$ d Y_t =
- \gamma \operatorname{sgn} (Y_t) dt + d \beta_t$$ où $\beta$ est un mouvement 
brownien
standard. \\ \\
Un tel processus admet une probabilité invariante, égale à la loi
d'une variable exponentielle symétrique de paramètre $2 \gamma$, et
son temps local en zéro, pris jusqu'à l'instant $t$, est
p.s. équivalent à $\gamma t$, quand $t$ tend vers l'infini. \\ \\ 
De plus, la propriété suivante nous sera utile par la suite : si $(\tilde{Y}_t)_{t \geq 0}$ est un
mouvement brownien avec drift $\gamma > 0$ issu de 0, et si on pose, pour tout
$t \in \mathbf{R}_+$, $S_t = \sup \{ \tilde{Y}_s , s \in [0,t] \}$, alors le
processus $(S_t - \tilde{Y}_t)_{t \geq 0}$ est la valeur absolue d'un
processus bang-bang de paramètre $\gamma$ (voir \cite{3}). \\ \\
Pour des discussions plus générales sur les processus de ce type, et en
particulier sur leur semi-groupe, voir également \cite{5}.       
\subsection{Définition des pénalisations étudiées et énoncé des
  théorèmes principaux de l'article} 
\noindent 
Après avoir défini la loi de l'araignée brownienne, nous lui
appliquons les changements de probabilité suivants : pour $\alpha
= (\alpha_i)_{i \in E}$ une famille de réels indexés par $E$, $\gamma
\in \mathbf{R}$ et $t \in \mathbf{R}_{+}$, on pose  
$$ ^{(\alpha,\gamma)}\mathbf{W}^{(t)} = \frac{\exp(\alpha_{N_t} X_t +
  \gamma L_t)}{\mathbf{W}_{(0,0)} [\exp(\alpha_{N_t} X_t +
  \gamma L_t)]} . \mathbf{W}_{(0,0)}$$
où $L_t$ est le temps local en 0 de $(X_s)_{s \leq t}$ : 
$$ L_t = \underset{\epsilon \rightarrow 0}{\lim \inf} \frac{1}{2
  \epsilon} \int_0^t \mathbf{1}_{X_s \leq \epsilon} ds$$ (en fait, la
  limite inférieure ci-dessus est presque sûrement une limite). \\ \\
Le but de notre article
  est de prouver les
  théorèmes suivants : \\ \\
\textbf{Théorème 1 : } \textit{ Il existe une mesure de probabilité
  $^{(\alpha,\gamma)}\mathbf{W}^{(\infty)}$ (sur la tribu
  $\mathcal{F}_{\infty}$ engendrée par les $\mathcal{F}_s$, $s \in
  \mathbf{R}_+$), telle que pour tout $s \in \mathbf{R}_+$ et tout $\Lambda_s
  \in \mathcal{F}_s$ : $$
  ^{(\alpha, \gamma)}   \mathbf{W}^{(t)}(\Lambda_s) \underset{t \rightarrow
  \infty}{\rightarrow} \, ^{(\alpha,\gamma)}\mathbf{W}^{(\infty)}
  (\Lambda_s)$$ De plus, on a :  
$$ ^{(\alpha,\gamma)}\mathbf{W}^{(\infty)}
  (\Lambda_s) = \mathbf{W}_{(0,0)} [\mathbf{1}_{\Lambda_s} M(\alpha,
  \gamma, X_s, N_s, L_s, s)]$$ où la fonction $M$ est donnée par le
  tableau suivant : }\\ \\ \\
\begin{tabular}{|p{3.5 cm}|c|}
\hline
Conditions sur $\alpha$, $\gamma$  & $M (\alpha,\gamma,x,k,l,s)$
\\
\hline
$\gamma \geq \alpha_m$ pour tout $m$ et $\gamma > 0$ & $e^{\gamma(l
  - x) - s \gamma^2/2}$ \\
\hline 
$\alpha_m = \max(\alpha) = \bar{\alpha}$ ssi $m \in J$ ($J \subset E$
et $J \neq \emptyset$), $\bar{\alpha} > \gamma$ et $\bar{\alpha} > 0$ &
$e^{\gamma l - s \bar{\alpha}^2/2} \left( e^{- \bar{\alpha} x} +
  \frac{\bar{\alpha} - \gamma}{\bar{\alpha} \underset{m \in J}{\sum}
    \mu_m} \sinh(\bar{\alpha} x) \mathbf{1}_{k \in J} \right)$ \\
\hline
$\gamma = 0$, $\alpha_m \leq 0$ pour tout $m \in E$ & 1 \\
\hline
$\alpha_m = 0$ si $m \in J$ ($J \subset E$ et $J \neq \emptyset$),
$\alpha_m < 0$ sinon, et $\gamma < 0$ & $e^{\gamma l} \left( 1 +
  \frac{|\gamma|}{\underset{m \in J}{\sum} \mu_m} x \mathbf{1}_{k
    \in J} \right)$ \\
\hline
$\alpha_m < 0$ pour tout $m \in E$ et $\gamma <0$ & $e^{\gamma l}
\left( 1 + \frac{\frac{1}{\alpha_{k}^2} + \underset{m \in E}{\sum}
    \frac{\mu_m}{\alpha_m \gamma}}{\underset{m \in E}{\sum} \mu_m
    \frac{|\alpha_m| + |\gamma|}{\alpha_m^2 \gamma^2}} x \right)$
\\
\hline
\end{tabular}
\\ \\ \\
\textit{En particulier, pour tout $s$, la restriction de
  $^{(\alpha, \gamma)} \mathbf{W}^{(\infty)}$ à $\mathcal{F}_s$
  est équivalente à la loi de l'araignée brownienne
  sur $[0,s]$, et la famille $(
  M(\alpha, \gamma, X_s, N_s, L_s, s))_{s \geq 0}$ des densités
  obtenues est une $\mathcal{F}_s$-martingale sous $\mathbf{W}_{(0,0)}$.}
  \\ \\ \\ 
\textbf{Théorème 2 : } \textit{Le processus canonique $(A_s)_{s \geq
  0}$ sous $^{(\alpha,\gamma)}\mathbf{W}^{(\infty)}$ peut
  être décrit de la manière suivante : \\ \\
- Si $\gamma > 0$ et $\gamma \geq \alpha_m$ pour tout $m$,
  $(X_s)_{s \geq 0}$ est la valeur absolue d'un processus bang-bang de
  paramètre $\gamma$, et la loi de $(N_s)_{s \geq 0}$
  conditionnellement à $(X_s)_{s \geq 0}$ est la même que sous
  $\mathbf{W}_{(0,0)}$ : les variables $(N_I)_{I \in \mathcal{I}}$
  ($\mathcal{I}$ étant l'ensemble des excursions de $X$) sont
  indépendantes de loi $\mu$. \\ \\
- Si $\bar{\alpha} = \max (\alpha) > \gamma$
  et $\bar{\alpha} > 0$, $(X_s)_{s \geq 0}$ est un processus dont
  la loi a une densité égale à $\frac{\bar{\alpha} -
    \gamma}{\bar{\alpha}} \exp(\gamma L_{\infty})$ par rapport à celle
  de la valeur absolue d'un mouvement brownien avec drift
  $\bar{\alpha}$ (dont $L_{\infty}$ est le temps local total sur tout
  $\mathbf{R}_+$), et $(N_s)_{s \geq 0}$ est obtenu en effectuant la même
  démarche que pour l'araignée initiale, puis en conditionnant le
  résultat par le fait que la dernière excursion de $(A_s)_{s
  \geq 0}$ se situe sur une branche $m$ vérifiant $\alpha_m =
  \bar{\alpha}$. \\ \\
- Si $\gamma = 0$ et $\alpha_m \leq 0$ pour tout $m$, 
  $(A_s)_{s \geq 0}$ est une araignée brownienne. \\ \\
- Si $\gamma < 0$ et $\alpha_m \leq 0$ pour tout $m$, on considère $(Y_s,R_s)_{s \geq 0}$ une araignée brownienne,
  $\mathbf{e}$ une variable exponentielle de paramètre $|\gamma|$ indépendante de
  $(Y_s,R_s)_{s \geq 0}$, $\tau_{\mathbf{e}}$ l'inverse du temps local
  de $(Y_s)_{s \geq 0}$ en $\mathbf{e}$, $(\hat{Y}_s)_{s \geq 0}$ un
  processus de Bessel de dimension 3 issu de 0 et indépendant des
  variables précédentes, $V$ une variable aléatoire (également
  indépendante des précédentes) définie sur $E$, et vérifiant les
  égalités suivantes pour $m \in E$ : 
  $$\mathbf{P} (V = m) = \frac{\mu_m}{\underset{k \in J}{\sum}
  \mu_k} \mathbf{1}_{m
  \in J}$$ si $J= \{ m \in E, \alpha_m = 0 \}$ est non vide, et
  $$\mathbf{P} (V = m) =
  \frac{ \mu_m \left( \frac{|\gamma|}{\alpha_m^2} + \underset{k \in
  E}{\sum} \frac{\mu_k}{|\alpha_k|} \right)} {\underset{k \in
  E}{\sum} \mu_k \frac{|\alpha_k| + |\gamma|}{\alpha_k^2}}$$ si $J =
  \emptyset$. \\ \\
Dans ces conditions, le processus $(X_s,N_s)_{s \geq 0}$ a même loi que
  $(\tilde{X}_s, \tilde{N}_s)_{s \geq 0}$, avec $(\tilde{X}_s, \tilde{N}_s) = (Y_s,R_s)$ pour $s \leq
  \tau_{\mathbf{e}}$, et  $(\tilde{X}_{s+ \tau_{\mathbf{e}}},
  \tilde{N}_{s+ \tau_{\mathbf{e}}}) = (\hat{Y}_{s}, V)$ pour $s \geq
  0$.} \\ \\ \\
Les Théorèmes 1 et 2 constituent une étude asymptotique complète des
  pénalisations exponentielles données au début de la section. \\
  On remarque que dans le cas où $\max \{\alpha_m, m \in E \} > 0$ et
  $\gamma = 0$, la densité de la restriction de 
$^{(\alpha,0)}\mathbf{W}^{(\infty)}$ à $\mathcal{F}_s$ ($s \geq
  0$), par rapport à celle de $\mathbf{W}_{(0,0)}$, est le produit
  d'une fonction de $s$ par une fonction de $A_s$. \\
Afin de comprendre ce résultat, on peut alors se demander pour quelles
  mesures de probabilités sur
  $\mathcal{C}_E$ une telle propriété a lieu. Le théorème suivant
  répond à la question
  dans le cas où les densités de probabilité sont suffisamment régulières. \\
  \\
\textbf{Théorème 3 : } \textit{ Soit $\nu$ une mesure de probabilité définie
  sur $\mathcal{C}_E$, différente de $\mathbf{W}_{(0,0)}$. \\
  On suppose que pour tout $s \geq 0$, la densité de la restriction de
  $\nu$ à $\mathcal{F}_s$ par rapport à celle de 
  $\mathbf{W}_{(0,0)}$ existe et
  s'écrit sous la forme : $$g(s,X_s,N_s) = h(s) f_{N_s}(X_s)$$
avec $f_m \in \mathcal{C}^2 (\mathbf{R}_+)$, $f_0(0)=f_m(0)=1$ pour tout $m
  \in E$, et $h \in \mathbf{C}^1(\mathbf{R}_+)$. \\
Dans ces conditions, il existe $\beta>0$ tel que $h(s) = e^{-s
  \beta^2/2}$ pour tout $s \geq 0$, et la mesure $\nu$ est une
  combinaison linéaire à coefficients positifs des mesures 
$^{(\alpha^{(m)},0)}\mathbf{W}^{(\infty)}$, où pour tout $m \in E$,
  $\alpha^{(m)}$ est donné par $\alpha^{(m)}_{m'} = \beta$ si $m=m'$ 
  et $\alpha^{(m)}_{m'} = 0$ sinon.}
  
\subsection{Interprétation heuristique des différents cas du Théorème
  2} 
\noindent
Les résultats donnés dans le Théorème 2 montrent que les processus
obtenus dépendent de manière assez complexe des paramètres $\alpha$ et
$\gamma$ définis précédemment. C'est pourquoi nous allons en donner
une interprétation heuristique. \\ \\
- Dans le premier cas du théorème, $\gamma > 0$ est le plus grand des paramètres de la 
pénalisation; de ce fait, son influence domine celle des $(\alpha_m)_{m 
\in E}$, et la loi limite obtenue ne dépend que de
$\gamma$. 
\\ Le processus canonique, sous cette loi limite, a 
alors tendance à rester près de l'origine, pour que son temps local en zéro
de $(X_t)_{t \geq 0}$ soit asymptotiquement plus grand. \\ Cette attraction vers l'origine correspond bien 
au comportement d'un processus bang-bang. \\ \\
- Dans le deuxième cas, l'influence qui domine est celle du plus grand
  coefficient $\bar{\alpha}$ : le processus canonique, sous la nouvelle loi de
  probabilité, reste (à partir d'un certain temps) dans une des
  branches $m \in E$ telles que $\alpha_m = \bar{\alpha}$. \\
De plus, la pénalisation exponentielle dominante est fortement 
liée à celle qui transforme un mouvement brownien standard en un 
mouvement brownien avec drift $\bar{\alpha}$, ce qui explique 
l'intervention de ce mouvement brownien avec drift dans la loi de 
$(X_t)_{t \geq 0}$.  \\ \\
- Dans le troisième cas, on pourrait penser que la pénalisation a tendance à 
empêcher le processus de trop s'éloigner de l'origine. \\ En réalité,
 la pénalisation étudiée est uniquement fonction de $(X_t,N_t)$, et le 
fait que l'on fasse tendre $t$ vers l'infini annule, à la limite, l'effet 
de cette pénalisation; le cas est analogue à celui d'un pont brownien sur 
$[0,t]$ ($t$ tendant vers l'infini) restreint à un intervalle fixé 
$[0,s]$ : ce processus tend, en loi, vers un mouvement brownien (voir 
\cite{9}). \\ \\
- Dans le dernier cas, la pénalisation du temps local à
  l'origine ($\gamma < 0$) domine, de sorte que le processus étudié
  reste dans une même branche à partir d'un certain temps; d'où 
l'intervention d'un processus de Bessel de dimension 3, qui n'est
autre qu'un mouvement brownien conditionné à rester positif sur tout 
$\mathbf{R}_+$. 
\subsection{Un petit guide de lecture de l'article}
\noindent
- Dans la suite de cet article, nous démontrons les trois théorèmes
 principaux, dans l'ordre où ils sont énoncés. \\ 
 Plus précisément, nous effectuons une étude préalable de la
  quantité : \\ $\mathbf{W}_{(x,k)} [\exp(\alpha_{N_t} X_t + \gamma
  L_t)]$ dans la Section 2, étude nécessaire à la preuve
  du Théorème 1 qui est achevée dans la Section 3. \\  Les Sections 4
  et 5 sont consacrées respectivement aux démonstrations des Théorèmes
  2 et 3. \\ \\
- On trouvera dans les preuves
  ci-dessous un certain nombre d'études de cas, selon les valeurs des
  différents paramètres. Une telle structure des démonstrations
  paraît inévitable, compte tenu du nombre assez important de ces
  paramètres. \\ Dans \cite{4} et \cite{9}, on peut également voir des situations où
  interviennent des distinctions de cas, analogues à celles
 rencontrées dans cet article. \\ \\
- Comme nous venons de l'évoquer ci-dessus, un certain nombre
d'estimations assez élémentaires (Propositions 2.1, 2.2, et 2.3,
Lemmes 2.4, 3.1 et 3.2), se ramenant assez rapidement à une
étude du mouvement brownien, sont faites préalablement aux
démonstrations des Théorèmes 1 et 2. \\
Nous conseillons au lecteur de faire une première lecture rapide de
ces estimations, puis de se concentrer sur les démonstrations des
théorèmes principaux de l'article, quitte à revenir ensuite sur la
preuve des résultats de la Section 2 et du Paragraphe 3.1. 
\section{Etude de l'expression $\mathbf{W}_{(x,k)} [\exp(\alpha_{N_t}
  X_t + \gamma L_t)]$}
\subsection{Enoncé des résultats obtenus}
\noindent
Afin de prouver l'existence de $^{(\alpha, \gamma)}\mathbf{W}^{(\infty)}$, nous allons commencer par définir une expression
qui majore $Z(\alpha, \gamma,x,k,t)=\mathbf{W}_{(x,k)} [\exp(\alpha_{N_t} X_t + \gamma
  L_t)]$ tout en étant équivalente à cette quantité quand $t$ tend
vers l'infini. \\ \\
Pour cela, introduisons les deux
quantités $I(\beta,
\gamma, x,t)$ et $J( \beta, x,t)$ ($\beta, \gamma \in \mathbf{R}$, $x,
t \in \mathbf{R}_+$), données par les égalités suivantes :
 $$ I(\beta, \gamma, x,t) = \mathbf{E}_x [\exp(\beta |Y_t| + \gamma
  L_t) \mathbf{1}_{T_0 \leq t}]$$
$$J(\beta, x,t) =
\mathbf{E}_x [\exp(\beta Y_t) \mathbf{1}_{T_0 > t}]$$
 où, sous $\mathbf{P}_x$, $(Y_t)_{t \geq
  0}$ est un mouvement brownien issu de $x$, $L_t$ le temps
local en zéro de $(Y_s)_{0 \leq s \leq t}$, et $T_0 = \inf \{s \geq 0,
Y_s = 0 \}$. \\ \\
De plus, posons : $$ J^*(\beta,x,t) = \sqrt{\frac{2}{\pi t^3}}
\frac{x}{\beta^2} \mathbf{1}_{\beta \neq 0} + \sqrt{\frac{2}{\pi t}} x
\mathbf{1}_{\beta=0} + 2 \sinh(\beta x) \exp(t \beta^2/2)
\mathbf{1}_{\beta > 0}$$ et définissons la quantité $I^*(\beta,
\gamma, x,t)$ par le tableau suivant : \\ \\ \\
\begin{tabular}{|c|c|} 
\hline
Conditions sur $\beta$ et $\gamma$ & $I^*(\beta, \gamma,x,t)$ \\
\hline
$\beta$, $\gamma < 0$ & $\sqrt{\frac{2}{\pi t^3}} \left(\frac{x}{\beta
    \gamma} + \frac{|\beta|+|\gamma|}{\beta^2 \gamma ^2} \right)$ \\
\hline
$\beta = 0$, $\gamma <0$ & $\frac{1}{|\gamma|} \sqrt{\frac{2}{\pi t}}$ \\
\hline
$\gamma = 0$, $\beta < 0$ & $\frac{1}{|\beta|} \sqrt{\frac{2}{\pi t}}$ \\
\hline 
$\beta = \gamma = 0$ & $1$ \\
\hline
$\beta > 0$, $\beta > \gamma$ & $\frac{1}{\beta - \gamma}
\sqrt{\frac{2}{\pi t}} + \frac{2 \beta}{\beta - \gamma} e^{-\beta x +
  t \beta^2/2}$ \\
\hline
$\gamma > 0$, $\gamma > \beta$ & $\frac{1}{\gamma - \beta}
\sqrt{\frac{2}{\pi t}} + \frac{2 \gamma}{\gamma - \beta} e^{-\gamma x
  + t \gamma^2/2}$ \\
\hline 
$\gamma = \beta > 0$ & $\gamma \sqrt{\frac{2 t}{\pi}} + 2(t \gamma^2 +
1) e^{- \gamma x + t \gamma^2/2}$ \\
\hline 
 
\end{tabular}
\\ \\ \\
Si on pose : $$ Z^*(\alpha,\gamma,x,k,t) = \underset{m \in E} {\sum} \mu_m I^*(\alpha_m, \gamma,x,t) +
J^*(\alpha_k,x,t)$$ on a alors les trois propositions suivantes : \\ \\
\textbf{Proposition 2.1 : } \textit{ Pour tous $\beta \in \mathbf{R}$ et $x \in
\mathbf{R}_+$ : \\ \\
$J(\beta, x, t) \leq J^*(\beta, x, t)$ pour tout $t \geq
  0$. \\ 
$J(\beta, x , t)$ est équivalent à $J^*(\beta,x,t)$ quand
$t$ tend vers l'infini.} \\ \\ \\
\textbf{Proposition 2.2 : } \textit{ Pour tous $\beta, \gamma \in \mathbf{R}$ et $x \in
\mathbf{R}_+$ : \\ \\
$I(\beta, \gamma,x, t) \leq I^*(\beta, \gamma, x, t)$ pour tout $t \geq
  0$. \\ 
$I(\beta, \gamma,x, t)$ est équivalent à $I^*(\beta, \gamma,x,t)$ quand
$t$ tend vers l'infini.} \\ \\ \\
\textbf{Proposition 2.3 : } \textit{Pour tous $\alpha \in \mathbf{R}^E$,
$\gamma \in \mathbf{R}$, $x \in \mathbf{R}_+$ et $k \in E$ : \\
\\
$Z(\alpha,\gamma,x,k,t) \leq Z^*(\alpha,\gamma,x,k,t)$ pour tout $t
\geq 0$. \\
$Z(\alpha,\gamma,x,k,t)$ est équivalent à $Z^*(\alpha,\gamma,x,k,t)$
quand $t$ tend vers l'infini.} \\ \\ \\
Remarquons tout de suite que les Propositions 2.1 et 2.2 entraînent la
Proposition 2.3. \\ 
En effet, on a : $$Z(\alpha,\gamma,x,k,t) = A_1 + A_2$$ avec 
$$ A_1 = \mathbf{W}_{(x,k)} [\exp(\alpha_{N_t} X_t + \gamma L_t)
\mathbf{1}_{T_0 \leq t}]$$
$$ A_2 = \mathbf{W}_{(x,k)} [\exp(\alpha_{N_t} X_t + \gamma L_t)
\mathbf{1}_{T_0 > t}]$$
où $T_0 = \inf \{ s \geq 0, X_s = 0 \}$. 
D'après la propriété d) de l'araignée (donnée au début de l'article), conditionnellement au
fait que $T_0 \leq t$, $N_t$ est une
variable de loi $\mu$, indépendante de $(X_t, L_t)$. Comme, d'autre
part, $(X_s)_{s \geq 0}$ sous $\mathbf{W}_{(x,k)}$ a même loi que $(|Y_s|)_{s \geq 0}$ sous
$\mathbf{P}_x$, on a : 
$$A_1 = \underset{m \in E}{\sum} \mu_m I(\alpha_m, \gamma, x,t)$$
Par ailleurs, si $(X_s)_{s \geq 0}$ ne s'annule pas avant $t$, il est
évident que $L_t = 0$ et $N_t = k$. \\
On a donc $A_2 = J(\alpha_k,x,t)$, et il en résulte l'égalité suivante : $$
Z(\alpha,\gamma,x,k,t) = \underset{m \in E} {\sum} \mu_m I(\alpha_m, \gamma,x,t) +
J(\alpha_k,x,t)$$
qui entraîne la Proposition 2.3, en supposant vraies les
Propositions 2.1 et 2.2. \\ \\
Il nous reste donc à démontrer ces deux propositions, ce qui est fait
dans les Paragraphes 2.2 et 2.3.
\subsection{Preuve de la Proposition 2.1} \noindent
Le principe de réflexion
 implique : \begin{align*}  J(\beta,x,t) &  = \mathbf{E}_x [e^{\beta Y_t}
 \mathbf{1}_{T_0 > t}] = 
\mathbf{E}_x [e^{\beta Y_t}
 \mathbf{1}_{Y_t > 0}] - 
 \mathbf{E}_x [e^{\beta Y_t}
 \mathbf{1}_{Y_t>0, T_0 \leq t}] \\  & = \mathbf{E}_x
 [e^{\beta Y_t} \mathbf{1}_{Y_t > 0}] - 
\mathbf{E}_x [e^{- \beta Y_t} \mathbf{1}_{Y_t < 0}] \\ & =
\frac{1}{\sqrt{2 \pi t}} \int_0^\infty (e^{- ((x-y)^2/2t) +
    \beta y} - 
e^{-((x+y)^2/2t) + \beta y}) dy \\ \end{align*} 
\textbf{Supposons $\beta <0$ : } De la majoration immédiate :
$$ e^{-(x-y)^2/2t} - e^{-(x+y)^2/2t} \leq \frac{(x+y)^2}{2t} -
\frac{(x-y)^2}{2t} = \frac{2xy}{t}$$
on déduit l'inégalité :
$$ J(\beta,x,t) \leq \sqrt{\frac{2}{\pi t^3}} x \int_0^\infty y
e^{\beta y} dy = \sqrt{\frac{2}{\pi t^3}} \frac{x}{\beta^2}=
J^*(\beta,x,t) \,.$$
Par ailleurs, on a les encadrements suivants : 
$$ 1 - \frac{(x-y)^2}{2t} \leq e^{-(x-y)^2/2t} \leq 1 -
\frac{(x-y)^2}{2t} + \frac{(x-y)^4}{8t^2}$$ 
$$ 1 - \frac{(x+y)^2}{2t} \leq e^{-(x+y)^2/2t} \leq 1 -
\frac{(x+y)^2}{2t} + \frac{(x+y)^4}{8t^2}$$
ce qui implique : $$  J^*(\beta,x,t)- J(\beta,x,t) \leq  \frac{1}{\sqrt{2 \pi t}} \int_0^\infty
\frac{(x+y)^4}{8 t^2} e^{\beta y} dy = x^5 t^{-5/2} C(\beta x)$$
où, pour tout $u<0$, $C(u)= \frac{1}{\sqrt{128 \pi}} \int_0^{\infty}
(1+y)^4 e^{uy} dy$ est fini. \\ \\ 
$J^*(\beta,x,t)$ est donc à la fois un majorant et un équivalent de
$J(\beta,x,t)$ quand $t \rightarrow \infty$ ($x$ étant fixé) : la
Proposition 2.1 est donc vraie pour $\beta <0$. \\ \\
\textbf{Supposons $\beta = 0$ : } On obtient ici $$ J(0,x,t) =
\frac{1}{\sqrt{2 \pi t}} \int_{-x}^x e^{-y^2/2t} dy$$ expression
admettant bien comme majorant et comme équivalent : $$ J^*(0,x,t) =
\sqrt{\frac{2}{\pi t}} x$$ quand $t$ tend vers l'infini.  \\ \\
\textbf{Supposons $\beta > 0$ : } On a l'égalité suivante : 
$$ \frac{1}{\sqrt{2 \pi t}} \int_{-\infty}^{\infty} (e^{
  -((x-y)^2/2t) + \beta y} -e^{-((x+y)^2/2t) + \beta y}) dy $$ $$ =  \frac{1}{\sqrt{2 \pi
  t}} \int_{-\infty}^{\infty} dz e^{-(z^2/2t) + \beta z} (e^{\beta x} -
  e^{- \beta x}) = 2 \sinh (\beta x) e^{t \beta ^2/2}$$ 
Or : $$  \frac{1}{\sqrt{2 \pi t}} \int_{-\infty}^{0} (e^{-((x-y)^2/2t) +
  \beta y}-e^{-((x+y)^2/2t) + \beta y}) dy = -J(- \beta,x,t)$$
D'où l'égalité : 
$$ J(\beta,x,t) = J(-\beta,x,t) + 2 \sinh(\beta x) e^{t \beta^2/2}$$
quantité qui admet comme majorant et comme équivalent : $$ J^*(\beta,x,t) = 2
\sinh(\beta x) e^{t \beta^2/2} + \sqrt{\frac{2}{\pi t^3}}
\frac{x}{\beta^2}$$ 
Nous venons donc de prouver la Proposition 2.1 dans tous les cas. \hfill$\Box$

\subsection{Preuve de la Proposition 2.2} \noindent
Afin de démontrer cette proposition, nous
allons donner quelques résultats sur la loi jointe de $(|Y_t|,L_t)$, lorsque $(Y_t)_{t \geq
  0}$ est un mouvement brownien issu de $x$ et $(L_t)_{t \geq 0}$ son
temps local en zéro. \\ 
Plus précisément, en notant (pour tout $x \in \mathbf{R}_+$), $\mathbf{P}_x$ la loi d'un mouvement
brownien réel issu de $x$, $(Y_t)_{t \geq 0}$ le processus canonique
de $\mathcal{C} (\mathbf{R}_+, \mathbf{R})$ et $L_t$ son temps local
en 0, nous allons prouver le lemme suivant : \\ \\
\textbf{Lemme 2.4 : } \textit{Avec les notations précédentes : \\ \\
- Pour tout $x \geq 0$, $\mathbf{P}_x [L_t+|Y_t| \in dz,
  L_t > 0] = \sqrt{\frac{2}{\pi t^3}} z(x+z) \exp \left( -
    \frac{(x+z)^2}{2t} \right) \mathbf{1}_{z > 0} dz$. \\ \\
- Conditionnellement au fait que $L_t> 0$, $\Theta_t =
\frac{|Y_t|}{L_t + |Y_t|}$ est une variable uniforme sur $[0,1]$,
indépendante de $L_t + |Y_t|$. \\ \\
Autrement dit, on a, pour $l, y > 0$ : $$ \mathbf{P}_x (L_t \in dl,
|Y_t| \in dy) = \sqrt{\frac{2}{\pi t^3}} (l+x+y) \exp \left( -
  \frac{(l + x+y)^2}{2t} \right) dy dl$$} \\
\noindent 
\textbf{Preuve : } En effectuant une intégration par rapport au premier et
au dernier temps d'annulation de $(Y_s)_{s \leq t}$, et en appliquant
la propriété de Markov au temps  $T_0 = \inf \{ t \geq 0, Y_t = 0 \}$,
on obtient, pour tous $y \in \mathbf{R}_+$ et $l>0$ :  
$$ \mathbf{P}_x (|Y_t| \in dy, L_t \in dl) = \underset{s_1+ s_2 \leq
  t}{\int} \mathbf{P}_x (T_0 \in ds_1) \mathbf{P}_0(|Y_{t-s_1}| \in
  dy) \, \, ...
  $$ $$... \, \,\mathbf{P}_0 (\underset{0 \leq u \leq t-s_1}{\sup} \{u | Y_u =
 0 \} \in d_{s_2}(t-s_2), L_{t-s_1} \in dl | |Y_{t-s_1}|=y)$$ 
Par un renversement du temps effectué sur le pont brownien :  
$$ \mathbf{P}_x (|Y_t| \in dy, L_t \in dl) $$ \begin{align*} & = \underset{s_1+s_2
  \leq t}{\int} \mathbf{P}_x (T_0 \in ds_1)  \mathbf{P}_y (T_0 \in ds_2,
  L_{t-s_1} \in dl, |Y_{t-s_1}| \in [0,dy])\\ & =  \underset{s_1+s_2
  \leq t}{\int} \mathbf{P}_x (T_0 \in ds_1) \mathbf{P}_y (T_0 \in ds_2)
  \mathbf{P}_0 (|Y_{t-s_1-s_2}| \in [0,dy], L_{t-s_1-s_2} \in dl) \\ 
& = \underset{s_1+s_2 \leq t}{\int} \frac{2 dy} {\sqrt{2 \pi
  (t-s_1-s_2)}}  \mathbf{P}_x (T_0 \in ds_1) \mathbf{P}_y (T_0 \in ds_2) 
 \mathbf{P}_0 (L_{t-s_1-s_2} \in dl | Y_{t-s_1-s_2} = 0) \end{align*}
Or la loi du temps local d'un pont brownien sur l'intervalle de temps $[0,t-s_1-s_2]$ est
  connue : c'est la loi (dite de Rayleigh) de la racine carrée d'une variable
  exponentielle de paramètre $\frac{1}{2(t-s_1-s_2)}$. \\ 
On en déduit : 
$$ \mathbf{P}_x (|Y_t| \in dy, L_t \in dl) =  \underset{s_1 + s_2
  \leq t}{\int} dy dl \frac{2 l e^{-l^2/2(t-s_1-s_2)}}{\sqrt{2 \pi
  (t-s_1-s_2)^3}} \mathbf{P}_x (T_0 \in ds_1) \mathbf{P}_y (T_0 \in ds_2)
  $$ \begin{align*} & = 2 dy dl \underset{s_1 + s_2 \leq t}{\int} \mathbf{P}_x(T_0 \in
  ds_1) \mathbf{P}_y(T_0 \in ds_2) D_l(t-s_1-s_2) \\ & = 2 dy dl \underset{s_1 + s_2
  \leq t}{\int} ds_1 ds_2 D_x(s_1) D_y(s_2) D_l(t-s_1-s_2) \\ & = 2
  dy dl D_x * D_y * D_l (t) = 2 D_{x+y+l} (t) dy dl \\ & =
  \sqrt{\frac{2}{\pi t^3}} (l+x+y) \exp \left( - \frac{(l+x+y)^2}{2t}
  \right) dy dl \end{align*}
$D_a(u)$ désignant la densité de $T_0$ en $u$ sous $\mathbf{P}_a$. \\ \\
Ces égalités impliquent le lemme annoncé (voir également \cite{6} pour une
  autre démonstration). \hfill$\Box$\\ \\
\textbf{Suite de la preuve de la Proposition 2.2} \\ \\
Avec les notations du Lemme 2.4, on peut écrire : \\ $\beta |Y_t| +
  \gamma L_t = (|Y_t| + L_t) (\gamma + (\beta - \gamma) \Theta_t)$, ce
  qui implique la formule suivante : $$ I(
  \beta, \gamma, x,t) = \mathbf{E} \left[ \sqrt{\frac{2}{\pi t ^3}}
  \int_0^{\infty} z(x+z) \exp \left( \frac{-(x+z)^2}{2t}+ z U_{\beta,\gamma} \right) dz \right]$$ où
  $U_{\beta, \gamma} = \gamma + (\beta-\gamma) \Theta_t$
 est une variable uniforme sur $[\beta, \gamma]$
  (ou bien $[\gamma, \beta]$ si $\gamma < \beta$). \\ 
Nous allons à présent distinguer plusieurs cas, selon les
  valeurs de $\beta$ et $\gamma$. \\ \\
\textbf{Supposons $\beta$, $\gamma < 0$ : } Dans ce cas, le théorème
  de convergence monotone prouve que $\mathbf{E} [\int_0^{\infty}
  z(x+z) e^{-((x+z)^2/2t) + z U_{\beta,\gamma} } dz]$ croît vers $\mathbf{E}
  [\int_0^{\infty} z(x+z) e^{z U_{\beta,\gamma} } dz ]$ quand $t$ tend vers
  l'infini. \\ \\ 
Or pour $\phi \in \mathbf{R}_-^*$, $\int_0^{\infty} z(x+z) e^{\phi z}
  dz = \frac{x} {\phi^2} + \frac{2}{|\phi|^3}$.  \\ \\
On en déduit que si $\beta \neq \gamma$ : \begin{align*} \mathbf{E} \left[
  \int_0^{\infty} z (x+z) e^{z U_{\beta,\gamma} } dz \right] & = \frac{1}{\gamma -
  \beta} \int_{\beta}^{\gamma} \left(\frac{x}{\phi^2} +
  \frac{2}{|\phi|^3} \right) d \phi \\ & = \frac{1}{\gamma - \beta}
  \left[-\frac{x}{\phi} + \frac{1}{\phi^2} \right]_{\beta}^{\gamma} =
  \frac{x}{\beta \gamma} + \frac{|\beta| +|\gamma|}{\beta^2
  \gamma^2} \end{align*} et que cette dernière égalité se prolonge en fait au cas
  où $\beta = \gamma$. \\ \\
Il en résulte que $I(\beta, \gamma, x,t)$ admet comme majorant et comme
  équivalent : $$ I^*( \beta, \gamma, x,t) = \sqrt{\frac{2}{\pi t^3}}
  \left(\frac{x}{\beta \gamma} + \frac{|\beta| + |\gamma|}{\beta^2
  \gamma^2} \right)$$  quand $t$ tend vers l'infini. \\ \\
\textbf{Supposons $\beta = 0$, $\gamma < 0$ : } On a, pour tout $z$ :
  $$\mathbf{E} [e^{z U_{\beta,\gamma}}] = \frac{1}{|\gamma|} \int_{\gamma}^0
  e^{\phi z} d \phi = \frac{1 - e^{\gamma z}}{|\gamma| z}$$ D'où :
  \begin{align*} I(\beta,
  \gamma, x,t) & = \frac{1}{|\gamma|} \sqrt{\frac{2}{\pi t^3}}
  \int_0^{\infty} (x+z) e^{-(x+z)^2/2t} dz \\ & - \frac{1}{|\gamma|}
  \sqrt{\frac{2}{\pi t^3}} \int_0^\infty (x+z) e^{- ((x+z)^2/2t) +
  \gamma z} dz \end{align*}
On a $\int_0^{\infty} (x+z) e^{\gamma z} dz < \infty$, donc le deuxième
  terme de l'expression ci-dessus, négatif, est dominé par $t^{-3/2}$
  quand $t$ tend vers l'infini. \\ \\
Par ailleurs, $$ \int_0^{\infty} (x+z) e^{-(x+z)^2/2t} dz = \left[ -t
  e^{-(x+z)^2/2t} \right]_0^{\infty} = t e^{-x^2/2t}$$ admet $t$ comme
  majorant et comme équivalent quand $t$ tend vers l'infini. \\ \\
Ces deux propriétés permettent d'en déduire que $I(\beta, \gamma,
  x,t)$ admet comme majorant et équivalent : $$ I^*(\beta, \gamma, x,t)
  = \frac{1}{|\gamma|} \sqrt{\frac{2}{\pi t}}$$
\textbf{Supposons $\gamma = 0$, $\beta < 0$ : } Par symétrie, ce cas est
  évidemment analogue au cas précédent. \\ \\
\textbf{Supposons $\beta=\gamma= 0$ : } On a : \begin{align*}  I(\beta, \gamma,
  x,t) & = \sqrt{\frac{2}{\pi t^3}} \int_0^{\infty} z(x+z)
  e^{-(x+z)^2/2t} dz \\ & = \sqrt{\frac{2}{\pi t^3}} \left[-tz
  e^{-(x+z)^2/2t} \right]_0^{\infty} + \sqrt{\frac{2}{\pi t}}
  \int_0^{\infty} e^{-(x+z)^2/2t} dz \\ & = 2 \mathbf{P}(\mathcal{N} \geq
  x/\sqrt{t}) \end{align*} où $\mathcal{N}$ est une variable gaussienne centrée
  réduite. \\ 
La Proposition 2.2 est donc vraie dans ce cas puisque $I^*(\beta, \gamma,
  x,t) = 1$. \\ \\
\textbf{Supposons $\beta > 0$ et $\beta > \gamma$ : } On a : $$
  \mathbf{E}[e^{z U_{\beta,\gamma} } ] = \frac{1}{\beta - \gamma}
  \int_{\gamma}^{\beta} e^{\phi z} d \phi = \frac{e^{\beta z} -
  e^{\gamma z}}{(\beta - \gamma) z}$$ 
On en déduit : \begin{align*}  I(\beta, \gamma, x,t) & = \frac{1}{\beta - \gamma}
  \sqrt{\frac{2}{\pi t^3}} \int_0^{\infty} (x+z) e^{-((x+z)^2/2t) +
  \beta z} dz \\ & - \frac{1}{\beta - \gamma} \sqrt{\frac{2}{\pi t^3}}
  \int_0^{\infty} (x+z) e^{-((x+z)^2/2t) + \gamma z} dz \end{align*}
Pour tout $\phi > 0$ : \begin{align*} &  \sqrt{\frac{2}{\pi t^3}} \int_0^{\infty}
  (x+z) e^{-((x+z)^2/2t) + \phi z} dz \\ = \, \,  & 
  \sqrt{\frac{2}{\pi t^3}} \left[
  -t e^{-((x+z)^2/2t) + \phi z} \right]_0^{\infty} + \phi
  \sqrt{\frac{2}{\pi t} } \int_0^{\infty} e^{-((x+z)^2/2t) + \phi z} dz
 \\ = \, \,  &  \sqrt{\frac{2}{\pi t}} e^{-x^2/2t} + 2 \phi e^{- \phi x + t
  \phi ^2/2} - \phi \sqrt{\frac{2}{\pi t}} \int_0^{\infty}
  e^{-((x-z)^2/2t) - \phi z} dz \end{align*} avec $$ \int_0^{\infty} e^{-((x-z)^2/2t)
  - \phi z} dz \leq \int_0^{\infty} e^{-\phi z} dz = \frac{1}{\phi}$$
Donc la quantité ci-dessus admet $ \sqrt{2/\pi t} + 2 \phi e^{-\phi x+
  t \phi^2/2}$ comme majorant et comme équivalent quand $t$ tend vers
  l'infini. \\
On peut en particulier en déduire que le second terme de $I(\beta,
  \gamma, x,t)$, négatif, est négligeable devant le premier quand $t$
  tend vers l'infini (quel que soit le signe de $\gamma$), ce qui
 permet de prendre : $$ I^*(\beta, \gamma, x, t) =
  \frac{1}{\beta- \gamma} \sqrt{\frac{2}{\pi t}} + \frac{2
  \beta}{\beta - \gamma} \exp(-\beta x + t \beta ^2/2)$$ comme
  majorant et équivalent de $I(\beta, \gamma,x,t)$. \\ \\ 
\textbf{Supposons $\gamma > 0$ et $\gamma > \beta$ : } Ce cas est
  analogue au précédent. \\ \\
\textbf{Supposons $\gamma = \beta > 0$ : } On a ici : $$ I(\beta,
  \gamma,x,t) = \sqrt{\frac{2}{\pi t^3}} \int_0^{\infty } z(x+z)
  e^{-((x+z)^2/2t)  + \gamma z} dz$$ 
Or $$\int_0^{\infty} (z(x+z) - \gamma t z -t) e^{-((x+z)^2/2t) + \gamma
  z} dz = \left[ -t z e^{-((x+z)^2/2t) + \gamma z} \right]_0^{\infty} =
  0$$
On en déduit : $$ I(\beta,\gamma,x,t) = \sqrt{\frac{2}{\pi t}}
  \int_0^{\infty} (\gamma z + 1) e^{-((x+z)^2/2t) + \gamma z} dz$$ 
Par ailleurs, on a : $$ \sqrt{\frac{2}{\pi t}} \int_0^{\infty}
  e^{-((x+z)^2/2t) + \gamma z} dz = 2 e^{- \gamma x + t \gamma^2/2} -
  \sqrt{\frac{2}{\pi t}} \int_0^{\infty} e^{-((x-z)^2/2t) - \gamma z} dz$$
quantité équivalente et inférieure à $2 e^{-\gamma x + t
  \gamma^2/2}$.  \\ \\
La quantité $- \gamma x \sqrt{\frac{2}{\pi t}} \int_0^{\infty}
  e^{-((x+z)^2/2t) + \gamma z} dz$ est donc négative et équivalente à
  $-2 \gamma x e^{- \gamma x + t \gamma^2/2}$. \\ \\
D'autre part, d'après un calcul précédemment effectué, \\ $\gamma
  \sqrt{\frac{2}{\pi t}} \int_0^{\infty} (x+z) e^{-((x+z)^2/2t) + \gamma
  z} dz$ est équivalent et inférieur à \\$\gamma \sqrt{\frac{2t}{\pi }}
  + 2 t \gamma^2 e^{-\gamma x + t \gamma ^2/2}$ (voir l'étude du cas
  $\beta > 0$ et $\beta > \gamma$). \\ \\
En additionnant les trois termes évalués ci-dessus, on obtient donc à
  nouveau :
  $$I^*(\beta, \gamma, x,t) = \gamma \sqrt{\frac{2t}{\pi}} + 2(t
  \gamma^2 + 1) \exp(-\gamma x + t \gamma^2/2)$$ comme majorant et
  équivalent pour $I(\beta, \gamma, x,t)$. \hfill$\Box$ \\ \\ \\
Nous venons donc d'achever la preuve des Propositions 2.1 et 2.2, qui
  entraînent la Proposition 2.3. 
\section{Preuve de l'existence de la mesure $^{(\alpha, \gamma)}\mathbf{W}^{(\infty)}$}
\subsection{Quelques lemmes techniques}
\noindent 
L'objet de la Section 3 est de prouver le Théorème 1. Pour cela, nous
aurons besoin de deux lemmes, dont le premier est le suivant : \\ \\
\textbf{Lemme 3.1 : } \textit{Pour tous $\beta$, $\gamma \in \mathbf{R}$, il existe
$D(\beta,\gamma)$ tel qu'on ait, pour tout $t \geq 1$ et tout $x \geq
0$ : $$ J^*(\beta,x,t) \leq D(\beta,\gamma) \sinh((\beta^{+} + 1) x)
I^*(\beta,\gamma,0,t)$$} \\
\noindent
\textbf{Preuve : } Fixons $\beta$ et $\gamma$ dans $\mathbf{R}$, $x$
dans $\mathbf{R}_+$, et supposons
$t \geq 1$. \\
- Si $\beta < 0$ et $\gamma < 0$, $J^*(\beta,x,t) = \sqrt{\frac{2}{\pi t
    ^3}} \frac{x}{\beta^2}$ et $ I^*(\beta, \gamma, 0,t) =
\sqrt{\frac{2}{\pi t^3}} \frac{|\beta|+|\gamma|}{\beta^2 \gamma^2}$,
ce qui implique : $$ J^*(\beta,x,t) = \frac{x \gamma^2}{|\beta| +
  |\gamma|} I^*(\beta, \gamma, 0,t)$$
- Si $\beta<0$ et $\gamma=0$, $ I^*(\beta,\gamma,0,t) = \frac{1}{|\beta|}
\sqrt{\frac{2}{\pi t}} \geq \frac{1}{|\beta|} \sqrt{\frac{2}{\pi t
^3}} $, et donc : $$ J^*(\beta,x,t) \leq \frac{x}{|\beta|}
I^*(\beta,\gamma,0,t)$$
- Si $\beta <0$ et $\gamma >0$, $I^*(\beta,\gamma,0,t) \geq
\frac{1}{\gamma - \beta} \sqrt{\frac{2}{\pi t}} \geq \frac{1}{\gamma -
  \beta} \sqrt{\frac{2}{\pi t^3}} $, d'où : $$ J^*(\beta,x,t) \leq x
\left(\frac{\gamma-\beta}{\beta^2} \right) I^*(\beta,\gamma,0,t)$$
- Si $\beta = 0$ et $\gamma < 0$, $J^*(\beta,x,t) = \sqrt{\frac{2}{\pi t}}
x$ et $I^*(\beta, \gamma,0,t) = \frac{1}{|\gamma|} \sqrt{\frac{2}{\pi
t}}$, d'où : $$ J^*(\beta, x,t) = |\gamma| x I^*(\beta, \gamma,0,t)$$
- Si $\beta = \gamma = 0$, $I^*(\beta,\gamma,0,t) = 1 \geq
\frac{1}{\sqrt{t}} \geq \sqrt{\frac{2}{\pi t}}$, et : $$ J^*(\beta, x,t)
\leq x I^*(\beta, \gamma, 0,t)$$ 
- Si $\beta = 0$ et $\gamma > 0$, $I^*(\beta,\gamma,0,t)  \geq
\frac{1}{\gamma-\beta} \sqrt{\frac{2}{\pi t}}$, d'où : $$ J^*(\beta,x,t)
\leq x (\gamma - \beta) I^*(\beta,\gamma,0,t)$$ 
- Si $\beta  > 0$ et $\gamma < \beta$, $J^*(\beta, x,t) =
\sqrt{\frac{2}{\pi t^3}} \frac{x}{\beta^2} + 2 \sinh (\beta x) e^{t
  \beta^2/2}$ et $I^*(\beta,\gamma,0,t) \geq \frac{1}{\beta-\gamma}
\sqrt{\frac{2}{\pi t^3}} + \frac{2 \beta}{\beta-\gamma} e^{t
  \beta^2/2}$. \\ \\
On en déduit que : $$ J^*(\beta,x,t) \leq \max
\left(\frac{x(\beta-\gamma)}{\beta^2},\sinh(\beta x)
  \frac{\beta-\gamma}{\beta} \right) I^*(\beta,\gamma,0,t)$$  
Or $x \leq \frac{\sinh(\beta x)}{\beta}$, d'où : $$ J^*(\beta,x,t) \leq
\max \left( \frac{\beta - \gamma}{\beta^3}, \frac{\beta -
    \gamma}{\beta} \right) \sinh(\beta x) I^*(\beta, \gamma, 0,t)$$
- Si $\beta > 0$ et $\gamma = \beta$, on a $I^*(\beta,\gamma,0,t) \geq
\beta \sqrt{\frac{2}{\pi t^3}} + 2 e^{t \beta^2/2}$. \\ \\ On obtient
donc : \begin{align*} J^*(\beta,x,t) & \leq \max \left(\frac{x}{\beta^3}, \sinh(\beta
  x) \right) I^*(\beta, \gamma,0,t) \\ & \leq \max \left(\frac{1}{\beta^4} ,1
\right) \sinh(\beta x) I^*(\beta,\gamma,0,t) \end{align*}
- Si $\beta > 0$ et $\gamma > \beta$, on a $I^*(\beta,\gamma,0,t) \geq
\frac{1}{\gamma-\beta} \sqrt{\frac{2}{\pi t^3}} + \frac{2
  \gamma}{\gamma - \beta}e^{t \beta^2/2}$ d'où : \begin{align*} J^*(\beta, \gamma,x)
 & \leq \max \left(x \frac{\gamma-\beta}{\beta^2}, \sinh(\beta x)
  \frac{\gamma-\beta}{\gamma} \right) I^*(\beta,\gamma,0,t) \\ & \leq \max
\left(\frac{\gamma - \beta}{\beta^3}, \frac{\gamma - \beta}{\gamma}
\right) \sinh(\beta x) I^*(\beta,\gamma,0,t)\end{align*}
Le Lemme 3.1 est donc prouvé dans tous les cas. \hfill$\Box$ \\ \\
Il est utilisé pour démontrer le lemme suivant : \\ \\
\textbf{Lemme 3.2 : } \textit{Pour tous $\alpha \in \mathbf{R}^E$, $\gamma \in
\mathbf{R}$, il existe $H(\alpha,\gamma)$, $\psi(\alpha) \in
\mathbf{R}^+$ tels que pour tous $x \in \mathbf{R}_+$, $k \in E$, et
$t \geq 1$ : $$Z^*(\alpha, \gamma, x,k,t) \leq H(\alpha,\gamma) \exp
(\psi (\alpha)x) Z^*(\alpha, \gamma,0,0,t)$$} \\
\noindent
\textbf{Preuve : }  
On observe, tout d'abord, que quels que soient $\beta$ et $\gamma$, il
existe $C(\beta, \gamma)$ tel que pour tous $x$, $t$ : $$I^*(\beta,
\gamma, x,t) \leq C(\beta, \gamma) (1+x) I^*(\beta, \gamma,0,t)$$ (en
fait, $I^*(\beta, \gamma, x, t) \leq I^*(\beta, \gamma, 0,t)$ dès que
$\sup(\beta, \gamma) \geq 0$). \\ \\
On en déduit que pour tous $\alpha = (\alpha_m)_{m \in E} \in
\mathbf{R}^E$, $\gamma \in \mathbf{R}$, $t, x \in \mathbf{R}_+$ :  \begin{align*} \underset{m \in E}{\sum} \mu_m I^*(\alpha_m, \gamma,x,t)
& \leq C (\alpha, \gamma) (1 + x) \underset{m \in E}{\sum} \mu_m
I^*(\alpha_m, \gamma, 0, t) \\ & \leq C(\alpha, \gamma) (1+x)
Z^*(\alpha,
\gamma, 0,0,t)  \tag{*}  \end{align*}
où $C(\alpha, \gamma) = \max \{C(\alpha_m, \gamma), m \in E \}$. \\ \\
A présent, posons $\delta(\alpha) = \max \{ \alpha_m^{+}, m \in E
\}$, $\nu = \min \{ \mu_m, m \in E \}$ ($\nu > 0$ puisque $\mu_m
> 0$ pour tout $m \in E$), et $D(\alpha, \gamma) = \max \{ D(\alpha_m,
\gamma), m \in E \}$. \\ \\
Si $t \geq 1$, le Lemme 3.1 permet d'obtenir les inégalités suivantes
: \begin{align*} J^*(\alpha_k,x,t) & \leq D(\alpha_k,
\gamma) \sinh((\alpha_k^+ + 1) x) I^*(\alpha_k,\gamma,0,t) \\ & \leq
D(\alpha,\gamma) \sinh ((\delta(\alpha) + 1) x) \underset{m \in
  E}{\sum} \frac{\mu_m}{\nu} I^*(\alpha_m, \gamma,0,t) \\ &
\leq \frac{D(\alpha, \gamma)}{\nu} \sinh((\delta(\alpha) +
1) x) Z^*(\alpha,\gamma, 0, 0,t) \end{align*} 
On en déduit, en utilisant l'inégalité (*) : $$ Z^* (\alpha, \gamma,x,k,t) \leq
\left(\frac{D(\alpha, \gamma)}{\nu} + C(\alpha,\gamma) \right)
\exp[(\delta(\alpha) + 1)x] Z^*(\alpha, \gamma,0,0,t) $$ ce qui achève
la démonstration du Lemme 3.2. \hfill$\Box$ \\ \\
Ce lemme nous permettra d'achever la preuve du Théorème 1, ce que nous
faisons dans le paragraphe suivant.
\subsection{Preuve du Théorème 1} 
\noindent
Soient $\alpha \in \mathbf{R}^E$, $\gamma \in \mathbf{R}$, $s \geq 0$
et $\Lambda_s \in \mathcal{F}_s$. En utilisant la propriété de Markov
au temps $s$, on obtient, pour tout $t \geq s$ : 
\begin{align*} ^{(\alpha,\gamma)}\mathbf{W}^{(t)} [\Lambda_s] & = \mathbf{W}_{(0,0)}
\left[ \mathbf{1}_{\Lambda_s} \frac{e^{\alpha_{N_t} X_t + \gamma
 L_t}}{\mathbf{W}_{(0,0)}[e^{\alpha_{N_t} X_t + \gamma L_t} ]}
\right] \\ & = \mathbf{W}_{(0,0)} \left[ \mathbf{1}_{\Lambda_s}
  \frac{\mathbf{W}_ {(0,0)}[e^{\alpha_{N_t} X_t + \gamma L_t} |
    \mathcal{F}_s ]}{\mathbf{W}_{(0,0)} [e^{\alpha_{N_t} X_t +
      \gamma L_t} ]} \right] \\ & = \mathbf{W}_ {(0,0)}\left[
  \mathbf{1}_{\Lambda_s} e^{\gamma L_s} \frac{\mathbf{W}_{(0,0)}
    [e^{\alpha_{N_t} X_t + \gamma(L_t - L_s)}|X_s,N_s
    ]}{\mathbf{W}_{(0,0)} [e^{\alpha_{N_t} X_t + \gamma L_t}]}
\right] \\ & = \mathbf{W}_ {(0,0)}\left[\mathbf{1}_{\Lambda_s}
  \exp(\gamma L_s) \frac{Z(\alpha, \gamma,X_s,N_s,t-s)}{Z(\alpha,
    \gamma,0,0,t)} \right] \end{align*} 
On sait que $\exp(\gamma L_s) \frac{Z (\alpha, \gamma, X_s,N_s,
  t-s)}{Z(\alpha, \gamma, 0,0,t)}$ est équivalent à $\exp(\gamma L_s)
\frac{Z^*(\alpha, \gamma,X_s,N_s,t-s)}{Z^*(\alpha, \gamma,0,0,t)}$
quand $t$ tend vers l'infini, $L_s$, $X_s$, $N_s$ étant fixés. \\ \\
Or $Z^*(\alpha, \gamma,x,k,u) = \underset{m \in E}{\sum} \mu_m
I^*(\alpha_m, \gamma,x,u) + J^*(\alpha_k,x,u)$ pour tous $x$, $k$, $u$,
donc d'après les expressions de $I^*$ et $J^*$ données precédemment, on a les
équivalents suivants : 
\\ \\ \\
\begin{tabular}{|p{3.5 cm}|c|} 
\hline 
Conditions sur $\alpha$, $\gamma$ & Equivalent de $Z^*(\alpha, \gamma,
x, k, u)$ quand $u \rightarrow \infty$ \\
\hline 
$\gamma \geq \alpha_m$ pour tout $m$, $\gamma = \alpha_m$ ssi $m \in
J$, $J$ sous-ensemble non vide de $E$, et $\gamma > 0$ & $ 2 \left(
  \underset{m \in J}{\sum} \mu_m \right) u \gamma^2 e^{-\gamma x + u
  \gamma^2/2}$ \\
\hline
$\gamma > \alpha_m$ pour tout $m \in E$ et $\gamma > 0$ & $\left( \underset{m
  \in E}{\sum} \frac{2 \gamma \mu_m}{\gamma- \alpha_m} \right) e^{-\gamma x
+ u \gamma^2/2}$ \\
\hline
$\alpha_m = \max(\alpha) = \bar{\alpha}$ ssi $m \in J$ ($J$
sous-ensemble non vide de $E$), $\bar{\alpha} > \gamma$ et
$\bar{\alpha} > 0$ & $e^{u \bar{\alpha}^2/2} \left( \frac{2
    \bar{\alpha}}{\bar{\alpha} - \gamma} \left( \underset{m \in
      J}{\sum} \mu_m \right) e^{- \bar{\alpha} x} + 2 \sinh
  (\bar{\alpha} x) \mathbf{1}_{k \in J} \right)$ \\
\hline
$\gamma = 0$, $\alpha_m = 0$ si $m \in J$ (sous-ensemble non vide de
$E$) et $\alpha_m < 0$ sinon & $\underset{m \in J}{\sum} \mu_m$ \\
\hline
$\gamma = 0$, $\alpha_m < 0$ pour tout $m \in E$ & $\sqrt{\frac{2}{\pi
    u}} \left( \underset{m \in E}{\sum} \frac{\mu_m}{|\alpha_m|}
\right)$ \\
\hline
$\alpha_m = 0$ si $m \in J$ (sous-ensemble non vide de $E$), $\alpha_m
< 0$ sinon, et $\gamma < 0$ & $\sqrt{\frac{2}{\pi u}} \left(
  \frac{\underset{m \in J}{\sum} \mu_m}{|\gamma|} + x \mathbf{1}_{k
    \in J} \right)$ \\
\hline
$\alpha_m < 0$ pour tout $m \in E$ et $\gamma < 0$ &
$\sqrt{\frac{2}{\pi u^3}} \left( \underset{m \in E}{\sum} \mu_m
  \frac{|\alpha_m| + |\gamma|}{\alpha_m^2 \gamma^2} + x \left(
    \frac{1}{\alpha_k^2} + \underset{m \in E}{\sum}
    \frac{\mu_m}{\alpha_m \gamma} \right) \right)$ \\
\hline
\end{tabular}
\\ \\ \\
On en déduit que l'expression $\exp(\gamma L_s) \frac{Z
  (\alpha, \gamma, X_s,N_s,t-s)}{Z(\alpha, \gamma,0,0,t)}$ converge,
quand $t$ tend vers l'infini, vers
$M(\alpha, \gamma, X_s, N_s, L_s,s)$. \\ \\ 
Par ailleurs, si $t \geq s + 1$, on a, d'après le Lemme 3.2, les
inégalités : \begin{align*}  Z(\alpha,
\gamma, X_s, N_s, t-s) & \leq Z^* (\alpha, \gamma, X_s, N_s, t-s) \\ &
\leq
H(\alpha, \gamma) e^{\psi(\alpha) X_s} Z^*(\alpha, \gamma,
0,0,t-s) \end{align*}
De plus : $$Z (\alpha, \gamma, 0,0,t) \geq \frac{1}{2} Z^*(\alpha,
\gamma,0,0,t)$$ pour $t$ assez grand, puisque $Z(\alpha, \gamma, 0,0,t)$ est équivalent à
$Z^*(\alpha, \gamma, 0,0,t)$ quand $t$ tend vers l'infini. \\ \\
D'autre part, pour $t$ assez grand : $$ \frac{Z^*(\alpha, \gamma,
  0,0,t-s)}{Z^*(\alpha, \gamma, 0,0,t)} \leq 2 M (\alpha,
\gamma,0,0,0,s) \leq 2$$
Il en résulte que pour $t$ assez grand : $$ e^{\gamma L_s}
\frac{Z(\alpha, \gamma, X_s, N_s,
  t-s)}{Z(\alpha, \gamma, 0,0,t)} \leq 4
H(\alpha, \gamma) \exp(\psi(\alpha) X_s + \gamma L_s)$$
Ce majorant étant intégrable sous $\mathbf{W}_{(0,0)}$, on en
déduit le Théorème 1, par convergence dominée. \hfill$\Box$
\section{Etude du processus associé à
  $^{(\alpha,\gamma)}\mathbf{W}^{(\infty)}$}
\noindent
Dans cette section, nous allons prouver le Théorème 2 en distinguant
plusieurs cas, selon l'expression de la
martingale $(M_s^{(\alpha,\gamma)})_{s \geq 0}$, \\ 
$M_s^{(\alpha,\gamma)}= M(\alpha, \gamma, X_s, N_s, L_s, s)$ étant la
densité de la restriction de $^{(\alpha,\gamma)}\mathbf{W}^{(\infty)}$ à
$\mathcal{F}_s$, par rapport à celle de $\mathbf{W}_{(0,0)}$ (nous
noterons plus simplement $M_s$ cette densité s'il n'y a pas
d'ambiguité possible).   
\subsection{Cas où $\gamma \geq \alpha_m$ pour tout $m$ et $\gamma > 0$}
\noindent
Sous $\mathbf{W}_{(0,0)}$, $(\tilde{Y}_s = L_s - X_s)_{s \geq 0}$ est un
mouvement brownien. La densité de la loi de $(\tilde{Y}_u)_{0 \leq u \leq s}$
sous $^{(\alpha,\gamma)}\mathbf{W}^{(\infty)}$, par rapport à celle
d'un mouvement brownien sur $[0,s]$, est donc égale à $M_s = \exp(\gamma
  \tilde{Y}_s - s \gamma^2/2)$. \\ \\
Par conséquent, $(\tilde{Y}_s)_{s \geq 0}$ est un mouvement brownien avec drift
$\gamma$ sous $^{(\alpha,\gamma)}\mathbf{W}^{(\infty)}$, et $(X_s)_{s
  \geq 0}$ peut être obtenu à partir de $(\tilde{Y}_s)_{s \geq
  0}$ grâce à l'expression : $X_s = \left( \underset{u \in
    [0,s]}{\sup} \tilde{Y}_u \right) - \tilde{Y}_s$. \\ \\
On en déduit que $(X_s)_{s \geq 0}$ a même loi que la valeur absolue d'un processus bang-bang de
paramètre $\gamma$. \\ \\ 
Par ailleurs, $N_s$ n'intervient pas dans l'expression de $M_s$, donc
le processus $(N_s)_{s \geq 0}$, sachant $(X_s)_{s \geq 0}$, a
même loi sous $^{(\alpha,\gamma)}\mathbf{W}^{(\infty)}$ que sous $\mathbf{W}_{(0,0)}$. 
\subsection{Cas où $\max \{ \alpha_m, m \in E \} > \max(\gamma,0)$ }  
\noindent
Dans ce paragraphe, nous posons $\bar{\alpha} = \max \{
\alpha_m, m \in E \}$ et nous notons $J$ l'ensemble (non vide) des éléments $m$ de $E$
tels que $\alpha_m = \bar{\alpha}$. Commen\c{c}ons tout d'abord par
étudier le cas particulier suivant : \\ \\
\textbf{a) Cas particulier : $\gamma = 0$ et il existe $m \in E$ tel
  que $J = \{ m \}$} \\ \\
On a dans ce cas : $$M_s = e^{-s \bar{\alpha}^2/2} \left( e^{- \bar{\alpha}
    X_s} + \frac{1}{\mu_m} \sinh(\bar{\alpha} X_s) \mathbf{1}_{N_s =
    m} \right)$$ pour tout $s \geq 0$. \\ \\
Considérons à présent un processus $(Y_t, R_t)_{t \geq 0}$ sur
$\mathbf{R}_E$ défini de la manière suivante : \\ \\
- $(Y_t)_{t \geq 0}$ est la valeur absolue d'un mouvement brownien
avec drift $\bar{\alpha}$. \\ \\
- Soit $\mathcal{I}$ l'ensemble des intervalles d'excursion de
$(Y_t)_{t \geq 0}$. Conditionnellement à $(Y_t)_{t \geq 0}$, $(R_t)_{t
  \geq 0}$ est constant sur
chaque intervalle $I \in \mathcal{I}$ ($R_t = R_I$ pour $t \in I$),
avec $R_I = m$ p.s. si $I$ est l'unique intervalle d'excursion non
borné de $(Y_t)_{t \geq 0}$, et avec les autres $(R_I)_{I \in
  \mathcal{I}}$ indépendants de loi $\mu$. \\ \\
Montrons alors que $^{(\alpha,0)}\mathbf{W}^{(\infty)}$ est la loi du
processus $(Y_t,R_t)_{t \geq 0}$. \\ \\
Pour cela, observons que si $t>0$, $k \in E$ et si $F$ est une
fonctionnelle mesurable bornée de $\mathcal{C}([0,t], \mathbf{R}_E)$
dans $\mathbf{R}$, on a : 
\begin{align*} \mathbf{E}
[F((Y_s)_{s \leq t}) \mathbf{1}_{R_t = k}] & = \mathbf{E} [F((Y_s)_{s
  \leq t}) \, \mathbf{P}
(R_t = k | (Y_s)_{s \in \mathbf{R}_+})] \\ & = \mathbf{E}
[F((Y_s)_{s \leq t}) \, (\mathbf{1}_{k=m} \mathbf{1}_{\forall s \geq
  t, Y_s > 0} + \mu_k \mathbf{1}_{\exists s \geq t, Y_s = 0})] \\ & =
\mathbf{E} [F((Y_s)_{s \leq t})  (\mathbf{1}_{k = m} \mathbf{P} (\forall s
\geq t, Y_s >0 | Y_t ) \\ & + \mu_k \mathbf{P} (\exists s \geq t, Y_s =
0 | Y_t ))] \end{align*} 
compte tenu de la propriété de Markov de $(Y_s)_{s \geq 0}$. \\ \\
Comme $Y_t= |B_t^{(\bar{\alpha})}|$ où
$(B_t^{(\bar{\alpha})})_{t \geq 0}$ est un mouvement brownien avec
drift $\bar{\alpha}$, on a : $$ \mathbf{P} (\exists s \geq t, Y_s
= 0 | Y_t)   = \mathbf{E} [ \mathbf{P} (\exists s \geq t, B_s^{(\bar{\alpha})} = 0
|B_t^{(\bar{\alpha})}) | | B_t^{(\bar{\alpha})}|] $$
Or $ \mathbf{P} (\exists s \geq t, B_s^{(\bar{\alpha})} = 0
|B_t^{(\bar{\alpha})})$ est égal à 1 si $B_t^{(\bar{\alpha})} \leq 0$
  et à $e^{-2 \bar{\alpha} B_t^{(\bar{\alpha})}}$ si
 $B_t^{(\bar{\alpha})} \geq 0$ : le premier cas est évident et le
  deuxième résulte du fait que  $(e^{-2 \bar{\alpha} B_{t \wedge
  T_0}^{(\bar{\alpha})}})_{t \geq s}$ est une martingale bornée si
  $T_0 = \inf \{ t \geq s, B_t^{(\bar{\alpha})} = 0 \}$. \\ \\
Compte tenu des densités en $Y_t$ et en $-Y_t$ de la loi de
  $B_t^{(\bar{\alpha})}$ : 
\begin{align*} 
   \mathbf{P} (\exists s \geq t, Y_s=0 |Y_t) 
  = \, \, &  \mathbf{P} (B_t^{(\bar{\alpha})} = -Y_t |Y_t) + e^{-2
  \bar{\alpha} Y_t}  \mathbf{P} (B_t^{(\bar{\alpha})} = Y_t |Y_t) \\
  = \, \, &  \frac{e^{- \bar{\alpha}
    Y_t}}{e^{\bar{\alpha} Y_t} + e^{- \bar{\alpha} Y_t}} +
  e^{-2 \bar{\alpha} Y_t}. \frac{e^{\bar{\alpha} Y_t}}{e^{\bar{\alpha} Y_t} +
  e^{- \bar{\alpha} Y_t}}  =
\frac{e^{-\bar{\alpha} Y_t}}{\cosh(\bar{\alpha} Y_t)} \end{align*} et
  donc $$ \mathbf{P}
(\forall s \geq t, Y_s > 0 | Y_t) = 1 - \frac{e^{- \bar{\alpha}
    Y_t}}{\cosh (\bar{\alpha} Y_t)} = \frac{\sinh(\bar{\alpha}
  Y_t)}{\cosh(\bar{\alpha} Y_t)}$$
On en déduit, en utilisant la densité de la loi de $(Y_s)_{s \leq t}$
  par rapport à celle du mouvement brownien réfléchi : \begin{align*}
 & \mathbf{E}[ F((Y_s)_{s \leq t}) \mathbf{1}_{R_t = k}] = 
  \mathbf{E} \left[F((Y_s)_{s \leq t}) \left(\mathbf{1}_{k=m} \frac{\sinh(\bar{\alpha}
  Y_t)}{\cosh(\bar{\alpha} Y_t)} + \mu_k \frac{e^{- \bar{\alpha}
    Y_t}}{\cosh(\bar{\alpha} Y_t)} \right) \right] \\ &  = 
  \mathbf{W}_{(0,0)} \left[ F((X_s)_{s \leq t}) \cosh (\bar{\alpha} X_t) e^{- t \bar{\alpha}^2/2} \left(
  \mathbf{1}_{k=m} \frac{\sinh(\bar{\alpha} X_t)}{\cosh(\bar{\alpha} X_t)}
  + \mu_k \frac{e^{-\bar{\alpha} X_t}}{\cosh(\bar{\alpha} X_t)}
  \right) \right] \\ &
  =     \mathbf{W}_{(0,0)} \left[ F((X_s)_{s \leq t})
  \mathbf{1}_{N_t=k} e^{- t \bar{\alpha}^2/2} \left( e^{- \bar{\alpha}
  X_t} + \frac{1}{\mu_m} \sinh(\bar{\alpha} X_t) \mathbf{1}_{N_t = m}
  \right) \right] \\
 & = \, \, ^{(\alpha, 0)} \mathbf{W}^{(\infty)} [F((X_s)_{s \leq t})
 \mathbf{1}_{N_t=k} ] \end{align*}
Nous avons donc montré que la loi de $((Y_s)_{s \leq t}, R_t)$ et
  celle de $((X_s)_{s \leq t}, N_t)$ sous $ ^{(\alpha, 0)}
  \mathbf{W}^{(\infty)}$ sont égales. \\ \\
Par ailleurs, la loi conditionnelle de $(R_s)_{s \leq t}$ sachant
  $((Y_s)_{s \leq t}, R_t)$ peut être décrite de la manière suivante :  
si $\mathcal{I}_t$ est l'ensemble des intervalles d'excursions de
  $(Y_s)_{s \leq t}$, $I_0$ l'élement de $\mathcal{I}_t$ contenant $t$,
  et pour tout $I \in \mathcal{I}_t$, $R_I=R_s$ avec $s$ quelconque dans
  $I$, alors $R_{I_0}= R_t$ p.s., et les variables
  $(R_I)_{I \in \mathcal{I}_t \backslash I_0}$ sont indépendantes de loi
  $\mu$. \\ \\
On déduit de cette description que la loi de $(R_s)_{s \leq t}$
  sachant $((Y_s)_{s \leq t}, R_t)$ est égale à celle de $(N_s)_{s
  \leq t}$ sachant $((X_s)_{s \leq t}, N_t)$ sous
  $\mathbf{W}_{(0,0)}$, et donc également sous
  $^{(\alpha,0)}\mathbf{W}^{(\infty)}$, puisque la densité de
  $^{(\alpha,0)}\mathbf{W}^{(\infty)}$ par rapport  à
  $\mathbf{W}_{(0,0)}$ ne dépend que de $X_t$ et
  $N_t$. \\ \\
Sous  $^{(\alpha,0)}\mathbf{W}^{(\infty)}$, on a donc d'une part
  l'égalité des lois de $((Y_s)_{s \leq t}, R_t)$ et de
 $((X_s)_{s \leq t}, N_t)$, d'autre part l'égalité des lois
  conditionnelles de $(R_s)_{s \leq t}$ sachant $((Y_s)_{s \leq t},
  R_t)$ et de  $(N_s)_{s \leq t}$ sachant $((X_s)_{s \leq t},
  N_t)$ : il en résulte l'égalité des lois de $(X_s,N_s)_{s \leq t}$
  et de $(Y_s,R_s)_{s \leq t}$. \\ \\
Le résultat annoncé au début du paragraphe est donc démontré,
 ce qui achève la démonstration du Théorème 2 dans le cas
particulier a). \\ \\
De plus, ce résultat entraîne les faits suivants, valables pour tout
  $s \geq 0$ sous
la probabilité $^{(\alpha,0)}\mathbf{W}^{(\infty)}$ : \\ \\
- Conditionnellement à $\mathcal{F}_s$ et au fait que $X_t$ ne
  s'annule pour aucun $t \geq s$, $L_{\infty} - L_s$ est nul. \\ 
- Conditionnellement à $\mathcal{F}_s$ et au fait que $X_t$ s'annule
  pour au moins un $t \geq s$, $L_{\infty} - L_s$ est une variable
  exponentielle de paramètre $\bar{\alpha}$. En particulier,
  $L_{\infty}$ est une variable exponentielle de paramètre $\bar{\alpha}$. \\ \\
Ces propriétés résultent du fait que le mouvement brownien avec drift
$\bar{\alpha}$ est un
processus fortement markovien dont le temps local total en zéro est une
variable exponentielle de paramètre $\bar{\alpha}$. \\ \\  
\noindent
\textbf{Remarque : } Si $E = \{ -1, 1 \}$, $\mu_1= \mu_{-1} = 1/2$ et $m=1$, le processus
$(X_t N_t)_{t \geq 0}$, qui est un mouvement brownien sous
$\mathbf{W}_{(0,0)}$, est un mouvement brownien avec drift
$\bar{\alpha}$ sous $^{(\alpha,0)}\mathbf{W}^{(\infty)}$. \\
Cela se vérifie aussi bien avec la martingale $(M_s)_{s \geq 0}$
qu'avec la description du processus $(Y_t,R_t)_{t \geq 0}$ donnée ci-dessus. \\ \\ 
Nous pouvons à présent traiter le cas suivant, plus général : \\ \\
\textbf{b) Cas où il existe $m \in E$ tel que $J = \{ m \}$, mais où
  $\gamma < \bar{\alpha}$ n'est pas nécessairement nul} \\ \\
On a maintenant : $$ M_s = \exp(\gamma L_s - s \bar{\alpha}^2/2)
\left( e^{- \bar{\alpha} X_s} + \frac{\bar{\alpha} -
    \gamma}{\bar{\alpha} \mu_m} \sinh( \bar{\alpha} X_s)
  \mathbf{1}_{N_s = m} \right)$$
D'autre part, sous $^{(\alpha,0)}\mathbf{W}^{(\infty)}$, $L_{\infty}$
est une variable exponentielle de paramètre
$\bar{\alpha}$; on peut donc définir la mesure de probabilité suivante :
 $$\nu = \frac{\bar{\alpha} - \gamma}{\bar{\alpha}} \exp(\gamma
  L_{\infty}). ^{(\alpha,0)}\mathbf{W}^{(\infty)}$$ (sous laquelle
$L_{\infty}$ est une variable exponentielle de paramètre $\bar{\alpha}
- \gamma$). \\ \\
Montrons que $\nu$ est exactement la mesure $^{(\alpha,\gamma)}\mathbf{W}^{(\infty)}$ que nous étudions, celle-ci étant donc absolument
continue par rapport à $^{(\alpha,0)}\mathbf{W}^{(\infty)}$. \\ \\
Pour prouver ce résultat, fixons $s \geq 0$ et $\Lambda_s \in
\mathcal{F}_s$. On a : 
\begin{align*} \nu(\Lambda_s) & = \,^{(\alpha,0)}\mathbf{W}^{(\infty)} \left[
  \frac{\bar{\alpha} - \gamma}{\bar{\alpha}} \mathbf{1}_{\Lambda_s}
  e^{\gamma L_s} e^{\gamma (L_{\infty} - L_s)} \right] \\ & = \,
  ^{(\alpha,0)}\mathbf{W}^{(\infty)} \left[ \frac{\bar{\alpha} -
  \gamma}{\bar{\alpha}} \mathbf{1}_{\Lambda_s} e^{\gamma L_s}
  \, ^{(\alpha,0)}\mathbf{W}^{(\infty)} [e^{\gamma (L_{\infty} -
  L_s)}|\mathcal{F}_s] \right] \\  = \mathbf{W}_{(0,0)} & \left[
  \frac{\bar{\alpha} - \gamma}{\bar{\alpha}} \mathbf{1}_{\Lambda_s}
  e^{\gamma L_s - s \bar{\alpha}^2/2} \left( e^{- \bar{\alpha} X_s} +
  \frac{1}{\mu_m} \sinh(\bar{\alpha} X_s) \mathbf{1}_{N_s =m} \right) \,
  ^{(\alpha,0)}\mathbf{W}^{(\infty)} [e^{\gamma (L_{\infty} - L_s)} |
  \mathcal{F}_s ] \right]  \end{align*}
\noindent
Par ailleurs, si $T_0 = \inf \{u \geq s, X_u = 0 \}$ et si $t \geq s$, on a :
  $$ ^{(\alpha,0)}\mathbf{W}^{(\infty)} [T_0 \leq t, \Lambda_s] =
  \mathbf{W}_{(0,0)} [M_t^{(\alpha,0)}. \mathbf{1}_{T_0 \leq t}
  . \mathbf{1}_{\Lambda_s}]$$ puisque $\{ T_0 \leq t \}$ et $\Lambda_s$ sont
  $\mathcal{F}_t$-mesurables. \\ \\
De plus, $ \{ T_0 \leq t \}$ et $\Lambda_s$ sont également $\mathcal{F}_{T_0
  \wedge t}$-mesurables, donc d'après le théorème d'arrêt :
  \begin{align*} ^{(\alpha,0)}\mathbf{W}^{(\infty)} [T_0 \leq t,
  \Lambda_s] & =
  \mathbf{W}_ {(0,0)} [M_{T_0 \wedge t}^{(\alpha, 0)} \mathbf{1}_{T_0
  \leq t} \mathbf{1}_{\Lambda_s}] \\ & = \mathbf{W}_{(0,0)}[M_{T_0}^{(\alpha, 0)}
  \mathbf{1}_{T_0 \leq t} \mathbf{1}_{\Lambda_s} ] \end{align*}
En faisant tendre $t$  vers l'infini, on a, par convergence monotone : \begin{align*}  ^{(\alpha,0)}\mathbf{W}^{(\infty)} [T_0 <
  \infty,\Lambda_s] &= \mathbf{W}_{(0,0)} [M_{T_0}^{(\alpha, 0)} \mathbf{1}_{T_0
  < \infty} \mathbf{1}_{\Lambda_s} ] \\ & = \mathbf{W}_{(0,0)} [M_{T_0}^{(
  \alpha, 0)} \mathbf{1}_{\Lambda_s}] \end{align*}  puisque $T_0 < \infty$ p.s. sous
  $\mathbf{W}_{(0,0)}$. \\ \\
Il en résulte : \begin{align*} ^{(\alpha,0)}\mathbf{W}^{(\infty)} [T_0<
  \infty, \Lambda_s] & =
  \mathbf{W}_{(0,0)} [\mathbf{1}_{\Lambda_s} \mathbf{W}_{(
  0,0)}[M_{T_0}^{( 
  \alpha,0)} | \mathcal{F}_s]] \\ & = \, ^{(\alpha,0)}\mathbf{W}^{(\infty)} \left[ \mathbf{1}_{\Lambda_s} \frac{\mathbf{W}_{(0,0)}
  [M_{T_0}^{(\alpha,0)} | \mathcal{F}_s]}{M_s^{(\alpha, 0)}}
  \right] \\ & = \, ^{(\alpha,0)}\mathbf{W}^{(\infty)} \left[
  \mathbf{1}_{\Lambda_s} \frac{\mathbf{W}_{(0,0)} [e^{-(T_0-s)
  \bar{\alpha}^2/2} | \mathcal{F}_s]}{e^{-\bar{\alpha} X_s} +
  \frac{1}{\mu_m} \sinh( \bar{\alpha} X_s) \mathbf{1}_{N_s = m}}
  \right] \end{align*}
Or, conditionnellement à $\mathcal{F}_s$, $T_0-s$ est le temps
  d'atteinte de zéro d'un mouvement brownien issu de $X_s$. On en déduit : $$ \mathbf{W}_{(0,0)}
  [e^{-(T_0-s) \bar{\alpha}^2/2} | \mathcal{F}_s] = e^{- \bar{\alpha}
  X_s}$$ et $$ ^{(\alpha,0)}\mathbf{W}^{(\infty)} [T_0< \infty, \Lambda_s] =
  \, ^{(\alpha,0)}\mathbf{W}^{(\infty)} \left[ \mathbf{1}_{\Lambda_s} \frac{e^{-
  \bar{\alpha} X_s}}{e^{- \bar{\alpha} X_s} + \frac{1}{\mu_m}
  \sinh(\bar{\alpha} X_s) \mathbf{1}_{N_s = m}} \right]$$ autrement
  dit : $$ ^{(\alpha,0)}\mathbf{W}^{(\infty)} [T_0<\infty |
  \mathcal{F}_s] = \frac{e^{-\bar{\alpha} X_s}}{e^{-\bar{\alpha} X_s}
  + \frac{1}{\mu_m} \sinh(\bar{\alpha} X_s) \mathbf{1}_{N_s = m}}$$
La loi conditionnelle de $L_{\infty}-L_s$, sachant la tribu engendrée
  par $\mathcal{F}_s$ et l'événement $\{T_0 < \infty\}$, a été décrite
 à la fin de l'étude du cas a). On déduit de cette description l'égalité
  suivante : $$ ^{(\alpha,0)}\mathbf{W}^{(\infty)}
  [e^{\gamma(L_{\infty}-L_s)} | \mathcal{F}_s] =
  \frac{\frac{\bar{\alpha}}{\bar{\alpha} - \gamma} e^{-\bar{\alpha}
  X_s} + \frac{1}{\mu_m} \sinh(\bar{\alpha} X_s) \mathbf{1}_{N_s =
  m}}{e^{- \bar{\alpha} X_s} + \frac{1}{\mu_m} \sinh(\bar{\alpha}
  X_s) \mathbf{1}_{N_s = m}}$$
\noindent
Il en résulte : \begin{align*}  \nu(\Lambda_s) & = \mathbf{W}_{(0,0)} \left[
  \mathbf{1}_{\Lambda_s} e^{\gamma L_s - s \bar{\alpha} ^2/2} \left( e^{-
  \bar{\alpha} X_s} + \frac{\bar{\alpha} - \gamma}{\bar{\alpha}
  \mu_m} \sinh(\bar{\alpha} X_s) \mathbf{1}_{N_s = m} \right) \right]
  \\ & = \mathbf{W}_ {(0,0)}[\mathbf{1}_{\Lambda_s} M_s^{(\alpha,
  \gamma)}] \end{align*}
On a donc l'égalité cherchée : $$\nu = \, ^{(\alpha,\gamma)}\mathbf{W}^{(\infty)}$$
qui implique le Théorème 2 dans le cas b). \\ \\
\textbf{c) Cas général} \\ \\
Une fois le Théorème 2 prouvé dans le cas particulier b), il est
  facile de l'étendre au cas général $\bar{\alpha} >
  \max(\gamma,0)$. \\ 
En effet, dans ce cas, la loi de
  $(X_t,N_t)_{t \geq 0}$ est une moyenne des
  lois données en b), avec une pondération
  $\frac{\mu_m}{\underset{k \in J}{\sum} \mu_k}$ pour chaque $m \in
  J$. \\ \\
On en déduit que le processus canonique sous
  $^{(\alpha,\gamma)}\mathbf{W}^{(\infty)}$ peut être décrit de la même
  manière qu'en b), sauf que sa
  dernière excursion se situe sur une branche quelconque appartenant à
  $J$, choisie aléatoirement à l'aide de la mesure $\mu$ : ceci
  correspond exactement à l'énoncé du Théorème 2. 
 \subsection{Cas où $\gamma < 0$ et $\alpha_m \leq 0$ pour tout $m \in E$}
  \noindent 
Dans ce cas, on a : $$ M_s = e^{\gamma L_s}(1 + \theta_{N_s} X_s)$$ où
  les $(\theta_k)_{k \in E}$, positifs, dépendant de $\alpha$, sont
  tels que : $$ \underset{k \in E}{\sum} \mu_k \theta_k = |\gamma|$$
Nous allons tout d'abord supposer qu'il existe $m \in E$ tel que $\theta_k =
  \frac{|\gamma|}{\mu_m}$ si $k=m$, et
  $\theta_k = 0$ si $k \neq m$. \\ \\
Dans ces conditions : $$M_s = e^{\gamma L_s} \left( 1 +
  \frac{|\gamma|}{\mu_m} X_s \mathbf{1}_{N_s=m} \right)$$ 
Considérons alors des réels positifs $l$ et $s$, une variable
  aléatoire $\mathcal{F}_{\tau_l}$-mesurable bornée $Y$ ($\tau_l$
  étant l'inverse, pris en $l$, du temps local de $(X_t)_{t \geq 0}$),
  et une fonctionnelle $F$ mesurable bornée de $\mathcal{C}([0,s],
  \mathbf{R}_E)$ dans $\mathbf{R}$. \\ \\
On a, lorsque $t \geq 0$ : \begin{align*} &  ^{(\alpha,\gamma)}\mathbf{W}^{(\infty)}
  \left[ \mathbf{1}_{\tau_l \leq t} Y F((X_{\tau_l + u}, N_{\tau_l +
  u})_{0 \leq u \leq s}) \right] \\ = \, \, & \mathbf{W}_{(0,0)} \left[
  M_{t+s}^{(\alpha, \gamma)} \mathbf{1}_{\tau_l \leq t} Y
  F((X_{\tau_l + u}, N_{\tau_l+ u})_{0 \leq u \leq s}) \right] \\ = \,
  \, &
  \mathbf{W}_{(0,0)} \left[ M_{(t+s) \wedge (\tau_l +
  s)}^{( 
  \alpha, \gamma)} \mathbf{1}_{\tau_l \leq t} Y F((X_{\tau_l + u},
  N_{\tau_l + u})_{0 \leq u \leq s}) \right] \\ = \, \, &
  \mathbf{W}_{(0,0)} \left[ M_{\tau_l + s}^{(\alpha, \gamma)}
  F((X_{\tau_l + u}, N_{\tau_l + u})_{0 \leq u \leq s} )
  \mathbf{1}_{\tau_l \leq t} Y \right] \end{align*}
en utilisant le théorème d'arrêt pour
  la deuxième égalité. \\ \\
Le théorème de convergence monotone entraîne alors, en faisant tendre $t$ vers
  l'infini : \begin{align*} 
 &  ^{(\alpha,\gamma)}\mathbf{W}^{(\infty)} \left[ \mathbf{1}_{\tau_l <
  \infty} Y F((X_{\tau_l + u}, N_{\tau_l + u})_{0 \leq u \leq s})
  \right] \\ = \, \, & e^{\gamma l} \mathbf{W}_{(0,0)} \left[e^{\gamma
  (L_{\tau_l + s} - L_{\tau_l})} \left( 1 + \frac{|\gamma|}{\mu_m}
  X_{\tau_l + s} \mathbf{1}_{N_{\tau_l + s} = m} \right) F((X_{\tau_l
  + u}, N_{\tau_l + u})_{0 \leq u \leq s}) Y \right] \end{align*}
compte tenu du fait que
  $\tau_l < \infty$ p.s. sous $\mathbf{W}_{(0,0)}$. \\ \\
D'après la propriété de Markov de l'araignée, $(X_{\tau_l + u},
  N_{\tau_l + u})_{0 \leq u \leq s}$ est indépendant de
  $\mathcal{F}_{\tau_l}$ sous $\mathbf{W}_{(0,0)}$ et a même loi que
  $(X_u, N_u)_{0 \leq u \leq s}$. \\ \\
On en déduit : $$ ^{(\alpha,\gamma)}\mathbf{W}^{(\infty)}
  \left[ \mathbf{1}_{L_{\infty} \geq l} Y F((X_{\tau_l +u}, N_{\tau_l
  + u})_{0 \leq u \leq s}) \right]$$ $$ =\exp(\gamma l) \,
  ^{(\alpha,\gamma)}\mathbf{W}^{(\infty)} [F(X_u,N_u)_{0 \leq u \leq
  s}] \mathbf{W}_{(0,0)} [Y]$$
En particulier, pour $F$ et $Y$ égaux à $1$, on obtient :
  $$^{(\alpha,\gamma)}\mathbf{W}^{(\infty)} [L_{\infty} \geq
  l] = \exp(\gamma l)$$
On a donc les caractéristiques suivantes : \\ \\
- $L_{\infty}$ est une variable exponentielle de paramètre
  $|\gamma|$. \\ \\
- Conditionnellement à $L_{\infty} \geq l$, $(X_s,N_s)_{0 \leq s \leq
  \tau_l}$ est une araignée brownienne arrêtée en $\tau_l$, et
  $(X_{\tau_l + s}, N_{\tau_l + s})_{s \geq 0}$ a pour loi
  $^{(\alpha,\gamma)}\mathbf{W}^{(\infty)}$; de plus, ces deux
  processus sont indépendants. \\ \\
On déduit de ce qui précède que conditionnellement à $L_{\infty} = l$,
  $(X_s,N_s)_{0 \leq s \leq \tau_l}$ est encore une araignée arrêtée
  en $\tau_l$, et $(X_{\tau_l + s},N_{\tau_l + s})_{s \geq 0}$ est un
  processus de loi $^{(\alpha,\gamma)}\mathbf{W}^{(\infty)}$,
  conditionné par le fait qu'il ne s'annule qu'au temps zéro, les deux
  processus étant encore indépendants. \\ \\
Pour déterminer la loi du deuxième processus, considérons $s \geq 0$, $\Lambda_s \in
  \mathcal{F}_s$, $l \geq 0$ et $t \geq s$. On a : \begin{align*}
  ^{(\alpha,\gamma)}\mathbf{W}^{(\infty)} [\Lambda_s, \tau_l \leq t]
  & =
  \mathbf{W}_{(0,0)} [M_t^{(\alpha, \gamma)} \mathbf{1}_{\Lambda_s}
  \mathbf{1}_{\tau_l \leq t}] \\ & = \mathbf{W}_{(0,0)} [M_{\tau_l \vee
  s}^{(\alpha, \gamma)} \mathbf{1}_{\Lambda_s} \mathbf{1}_{\tau_l \leq
  t}] \end{align*} d'où $$^{(\alpha,\gamma)}\mathbf{W}^{(\infty)} [\Lambda_s, \tau_l
  < \infty] = \mathbf{W}_ {(0,0)}[M_{\tau_l \vee s}^{(\alpha,
  \gamma)} \mathbf{1}_{\Lambda_s}]$$ et donc : \begin{align*} & ^{(\alpha,\gamma)}\mathbf{W}^{(\infty)} [\Lambda_s, L_{\infty} \leq l] = \mathbf{W}_{(0,0)}
  [(M_s^{( \alpha, \gamma)} - M_{\tau_l \vee s}^{(\alpha,
  \gamma)}) \mathbf{1}_{\Lambda_s} ] \\ & = \mathbf{W}_{(0,0)} \left[
  \left( e^{\gamma L_s} \left( 1 + \frac{|\gamma|}{\mu_m} X_s
  \mathbf{1}_{N_s = m} \right) - e^{\gamma l} \right) \mathbf{1}_{L_s
  \leq l} \mathbf{1}_{\Lambda_s} \right] \\ & = \mathbf{W}_{(0,0)}
  [L_s \leq l] \mathbf{W}_ {(0,0)}\left[ \mathbf{1}_{\Lambda_s} \left(
  e^{\gamma L_s} \left( 1 + \frac{|\gamma|}{\mu_m} X_s
  \mathbf{1}_{N_s = m} \right) - e^{\gamma l} \right) | L_s \leq l
  \right] \end{align*}
Comme $^{(\alpha,\gamma)}\mathbf{W}^{(\infty)} [L_{\infty} \leq l] = 1
  - e^{\gamma l}$, on a : \begin{align*} & ^{(\alpha,\gamma)}\mathbf{W}^{(\infty)}
  [\Lambda_s | L_{\infty} \leq l] \\ & = \frac{\mathbf{W}_{( 0,0)} [L_s \leq
  l]}{1-e^{\gamma l}} \mathbf{W}_{(0,0)}\left[ \mathbf{1}_{\Lambda_s}
  \left( e^{\gamma L_s} \left( 1 + \frac{|\gamma|}{\mu_m} X_s
  \mathbf{1}_{N_s = m} \right) - e^{\gamma l} \right) | L_s \leq l
  \right] \\ & = \frac{\mathbf{W}_ {(0,0)}[L_s \leq l]}{1 -
  e^{\gamma l}} \tilde{\mathbf{W}}_s (l) \left[ \mathbf{1}_{\Lambda_s}
  \left( e^{\gamma L_s} \left( 1 + \frac{|\gamma|}{\mu_m} X_s
  \mathbf{1}_{N_s = m} \right) - e^{\gamma l} \right) \right] \end{align*} où
  $\tilde{\mathbf{W}}_s (l)$ est la loi de $(X_u,N_u)_{u \leq s}$
  conditionné par l'événement $\{ L_s \leq l \}$. \\ \\
Quand $l$ tend vers zéro, $\frac{\mathbf{W}_{(0,0)} [L_s \leq l]}{1
  - e^{\gamma l}}$ tend vers $\frac{1}{|\gamma|} \sqrt{\frac{2}{\pi
  s}}$. \\ \\
D'autre part, si $L_s \leq l$ et si $(X_s,N_s)$ est fixé, 
  $e^{\gamma L_s} \left( 1 + \frac{|\gamma|}{\mu_m} X_s
  \mathbf{1}_{N_s = m} \right) - e^{\gamma l}$ tend vers 
  $\frac{|\gamma|}{\mu_m} X_s \mathbf{1}_{N_s = m}$ quand $l$ tend vers
  zéro. \\ \\
Ceci permet de démontrer : $$ ^{(\alpha,\gamma)}\mathbf{W}^{(\infty)}
  [\Lambda_s | L_{\infty} = 0] = \tilde{\mathbf{W}}_s (0)
  \left[\mathbf{1}_{\Lambda_s} \sqrt{\frac{2}{\pi s}} \frac{X_s}{\mu_m}
  \mathbf{1}_{N_s = m} \right]$$ où $\tilde{\mathbf{W}}_s (0)$ est la
  loi d'une araignée sur $[0,s]$, conditionnée par sa non-annulation
  en dehors du temps $0$; sous $\tilde{\mathbf{W}}_s
  (0)$, $(X_u)_{u \leq s}$ est un méandre brownien de durée $s$. \\ \\
On en déduit que sous $^{(\alpha,\gamma)}\mathbf{W}^{(\infty)}$,
  et conditionnellement au fait que $X_s > 0$ pour tout $s > 0$,
  $(X_s)_{s \geq 0}$ est un processus de Bessel de dimension 3, et
  $N_s =m$ pour tout $s$. \\ \\
On a donc la description de $(X_t,N_t)_{t \geq 0}$ dans le cas où un
  seul des $\theta_k$ précédemment donnés est nul. \\ \\
Le cas général est simple à étudier à présent; en effet, il suffit de
  faire une moyenne pondérée des mesures précédemment décrites pour
  chacun des $m \in E$ (avec la pondération $\frac{\mu_m
  \theta_m}{|\gamma|}$). \\ \\
\textbf{Remarque : } Dans le cas étudié ci-dessus ($\gamma <0$ et
  $\alpha_m \leq 0$ pour tout $m$), on peut
  vérifier directement, à partir
  de l'expression de $M_s$, que la loi du processus $(X_u)_{u \leq s}$ a
  une densité $e^{\gamma L_s} (1 + |\gamma| X_s)$ par rapport à la
  loi d'un mouvement brownien réfléchi sur $[0,s]$. \\ 
La validité de la description de $(X_s)_{s \geq 0}$ donnée dans le
  Théorème 2 peut alors se déduire du Théorème 1.1. de \cite{9} (appliqué à
  la fonction $\phi(y) = |\gamma| e^{\gamma y}$) et du théorème
  d'équivalence de Lévy. \\ 
Par ailleurs, dans le cas où $\theta_k > 0$ ssi $k=m$, on a :
 $$^{(\alpha,\gamma)} \mathbf{W}^{(\infty)} [N_t \neq m]
  = \mathbf{W}_{(0,0)} [e^{\gamma L_t} \mathbf{1}_{N_t \neq m} ]
  \underset{t \rightarrow \infty}{\rightarrow} 0$$
ce qui donne une preuve rapide du fait que l'excursion non bornée de
  $(A_s)_{s \geq 0}$ se situe presque sûrement sur la demi-droite d'indice $m$. 
\subsection{Cas où $\gamma = 0$ et $\alpha_m \leq 0$ pour tout $m \in E$}
  \noindent
Ce cas est le plus simple de tous : $^{(\alpha,0)}\mathbf{W}^{(\infty)}$
  est exactement la loi d'une araignée brownienne, puisque $M_s^{(
  \alpha, 0)}$ est constante et égale à 1. \\ \\ \\
A présent, nous venons d'achever la preuve du Théorème 2 dans tous les
  cas possibles. \hfill$\Box$ \\ \\ \\
\textbf{Remarque : } Le Théorème 2
  indique différents comportements possibles pour le processus limite,
 selon les valeurs des réels $\alpha_m$ ($m \in E$) et
  $\gamma$. \\
Cette distinction de cas généralise celle que B. Roynette, P. Vallois
  et M. Yor obtiennent dans \cite{9} lors de l'étude des pénalisations
  exponentielles du mouvement brownien. Une distinction de cas du même
  type apparaît également dans les résultats prouvés par Y. Hariya et
  M. Yor dans \cite{4}. \\ \\
Par ailleurs, il pourrait être intéressant d'étudier d'autres
  pénalisations de l'araignée brownienne, liées par exemple aux temps
  passés par l'araignée dans les différentes branches, dont la loi
  jointe, sur un intervalle de temps fixe, est décrite en \cite{2}. \\ \\
\section{Preuve du Théorème 3} 
\noindent
 Soit $\nu$ une mesure de probabilité vérifiant les conditions de
 l'énoncé du Théorème 3. \\
La famille des variable aléatoires $(g(s,X_s,N_s))_{s \geq 0}$ est une 
martingale sous
 $\mathbf{W}_{(0,0)}$, ce qui implique les résultats suivants, compte
 tenu du semi-groupe de l'araignée brownienne : \\ \\
- Pour tous $x \in \mathbf{R}_+^*$, $s \geq 0$ et $m \in E$,
  $\frac{\partial g}{\partial s} (s,x,m) + \frac{1}{2}
  \frac{\partial^2 g}{\partial x^2} (s,x,m) = 0$. \\ \\
- Pour tout $s \geq 0$, $\underset{m \in E}{\sum} \mu_m
  \frac{\partial g}{\partial x} (s,0,m) = 0$. \\ \\
La première égalité donne : $$ h'(s)
  f_m(x) + \frac{1}{2} h(s) f''_m(x) = 0$$
Comme $h(s)$ est non nul pour tout $s$, on en déduit : 
$$ f''_m(x) = - \frac{2 h'(s)}{h(s)} f_m(x)$$
ce qui implique que $C=\frac{h'(s)}{h(s)}$ ne dépend pas de $s$. \\ \\
Si on suppose $C > 0$, $f''_m(x)$ est négatif pour tout $x$, et $f'_m$ est
décroissante. Comme $f_m$ est une fonction positive, la limite de
$f'_m(x)$ quand $x$ tend vers l'infini est positive, et $f_m$ est
croissante : pour tout $x$, $f_m(x) \geq f_m(0) =1$. \\
On en déduit que $f''_m (x) = -2C f_m(x) \leq - 2C$, ce
qui est contradictoire avec la positivité de $f_m$. \\ \\
Si on suppose $C=0$, toutes les fonctions $f_m$ ($m \in E$) sont
affines et positives : $f_m(x) = 1 + \lambda_m x$ où $\lambda_m \geq
0$. \\ Or, pour tout $s \geq 0$, $\underset{m \in E}{\sum} \mu_m
  \frac{\partial g}{\partial x} (s,0,m) = 0$, ce qui implique
  $\underset{m \in E}{\sum} \lambda_m \mu_m = 0$, d'où $\lambda_m = 0$
  pour tout $m$. \\
La mesure $\nu$ est alors égale à
$\mathbf{W}_{(0,0)}$, ce qui est exclu dans l'hypothèse du Théorème
3. \\ \\
Nous venons donc de prouver que $C$ est strictement
négatif, notons le $- \beta^2/2$, avec $\beta > 0$. \\ 
Comme $h(0)=1$ (puisque $f_0(0)=1$), on a $h(s)= e^{-s \beta^2/2}$ comme
annoncé, et $f''_m(x) = \beta^2 f_m (x)$. \\ \\
On en déduit qu'il existe $\delta_m$ et $\lambda_m \in \mathbf{R}$
  tels que : $$f_m(x) = \delta_m \exp(-\beta x) + \lambda_m \sinh (\beta
  x)$$ pour tout $x \geq 0$; comme $f_m(0)=1$, $\delta_m=1$ pour tout
  $m$. \\ \\
Par ailleurs, $f_m(x) \geq 0$ pour tout $x \geq 0$, donc $\lambda_m
  \geq 0$ pour tout $m$. \\ \\
De plus, on doit avoir, pour tout $s \geq 0$, $\underset{m \in E}{\sum} \mu_m
  \frac{\partial g}{\partial x} (s,0,m) = 0$, ce qui implique $\underset{m \in
  E}{\sum} \mu_m f'_m(0) = 0$, soit $\underset{m \in E}{\sum} \mu_m
 (1 - \lambda_m) = 0$ et $\underset{m \in E}{\sum} \mu_m \lambda_m =
  1$. \\ \\
On en déduit que $\nu$ est une moyenne pondérée des mesures
  $^{(\alpha^{(m)},0)}\mathbf{W}^{(\infty)}$, la pondération étant 
$\mu_m \lambda_m$, ce qui achève la preuve du Théorème 3. \hfill$\Box$
\\ \\
Les processus associés aux mesures vérifiant l'énoncé
du Théorème 3 peuvent être considérés comme des généralisations du
 mouvement brownien avec drift.  

\bibliographystyle{alpha}
 \bibliography{araignee3}

\end{document}